\documentclass[11pt]{amsart}

\newtheorem*{question}{Question}
\newtheorem{theorem}{Theorem}

\newtheorem{corollary}[theorem]{Corollary}
\newtheorem{lemma}[theorem]{Lemma}
\newtheorem{ex}[theorem]{Example}
\newtheorem{rem}[theorem]{Remark}

\newcommand{\ci}{\circle*{.15}}
\renewcommand{\mp}{\multiput}
\newcommand{\num}[1]{{\raisebox{-.5\unitlength}{\makebox(0,0)[b]{${#1}$}}}}

\newcommand{\loops}{\mathcal{L}_{G}}
\newcommand{\quot}{\widetilde{W}/W}
\newcommand{\core}{\mathfrak{c}}
\newcommand{\bounded}{\mathfrak{b}}
\newcommand{\boundedtominreps}{\mathfrak{a}}
\newcommand{\minreps}{\widetilde{W}^{S}}
\newcommand{\Bn}{B_{n}}
\newcommand{\Cn}{C_{n}}
\newcommand{\Dn}{D_{n}}
\newcommand{\id}{\mathrm{id}}
\newcommand{\segment}[2]{\Sigma_{#2}^{#1}}
\newcommand{\fragment}[1]{F^{#1}}

\newcommand{\segmentD}[2]{\Sigma_{#2}^{#1}}
\newcommand{\partitions}{\mathcal{P}}
\newcommand{\given}{\, | \,} 
\newcommand{\tchoose}[2]{
\left[\begin{smallmatrix} #1 \\#2 
        \end{smallmatrix} \right]_{t}}
\newcommand{\chs}[2]{
\left(\begin{smallmatrix} #1 \\#2 
        \end{smallmatrix} \right)}

\newlength{\cellsize} 
\newsavebox{\cell}
\sbox{\cell}{\begin{picture}(18,18)
\put(0,0){\line(1,0){18}}
\put(0,0){\line(0,1){18}}
\put(18,0){\line(0,1){18}}
\put(0,18){\line(1,0){18}}
\end{picture}}
\newcommand\cellify[1]{\def\thearg{#1}\def\nothing{}%
\ifx\thearg\nothing
\vrule width0pt height\cellsize depth0pt\else
\hbox to 0pt{\usebox{\cell} \hss}\fi%
\vbox to \cellsize{
\vss
\hbox to \cellsize{\hss$#1$\hss}
\vss}}
\newcommand\tableau[1]{\vtop{\let\\\cr
\baselineskip -16000pt \lineskiplimit 16000pt \lineskip 0pt
\ialign{&\cellify{##}\cr#1\crcr}}}
\begin{document}
\title{Affine partitions and  affine Grassmannians} 
\date{\today}
\author{Sara C. Billey \and Stephen A. Mitchell}

\address{Department of Mathematics, University of Washington, Seattle, WA}
\email{billey@math.washington.edu}
\email{mitchell@math.washington.edu}
\thanks{S.B. was supported by UW Royalty Research Grant.}
\thanks{S.M. was supported by the National Science Foundation.}

\begin{abstract}
We give a bijection between certain colored partitions and the
elements in the quotient of an affine Weyl group modulo its Weyl
group.  By Bott's formula these colored partitions give rise to some
partition identities. In certain types, these identities have
previously appeared in the work of Bousquet-Melou-Eriksson,
Eriksson-Eriksson and Reiner. In other types the identities appear to
be new.  For type $A_{n}$, the affine colored partitions form another
family of combinatorial objects in bijection with $n+1$-core
partitions and $n$-bounded partitions. Our main application is to
characterize the rationally smooth Schubert varieties in the affine
Grassmannians in terms of affine partitions and a generalization of
Young's lattice which refines weak order and is a subposet of Bruhat
order. Several of the proofs are computer assisted.
\end{abstract}

\maketitle

\section{Introduction} \label{s:intro}

Let $W$ be a finite irreducible Weyl group associated to a simple
connected compact Lie Group $G$, and let $\widetilde{W}$ be its
associated affine Weyl group.  In analogy with the Grassmannian
manifolds in classical type $A$, the quotient $\quot$ is the indexing
set for the Schubert varieties in the affine Grassmannians $\loops$.
Let $\minreps$ be the minimal length coset representatives for
$\quot$.  Much of the geometry and topology for the affine
Grassmannians can be studied from the combinatorics of $\minreps$ and
vice versa.  For example, Bott \cite{Bott-56} showed that the
Poincar\'e series for the cohomology ring for the affine Grassmannian
is equal to the length generating function for $\minreps$, and this
series can be written in terms of the \textit{exponents}
$e_{1},e_{2},\dotsc, e_{n}$ for $W$ as
\begin{equation}\label{e:bott.formula}
P_{\minreps} (t) = \frac{1}{(1-t^{e_{1}}) (1-t^{e_{2}}) \cdots
(1-t^{e_{n}})}.
\end{equation}

Bott's formula suggests there is a natural bijection between elements
in $\minreps$ and a subset of partitions that preserves length.  The
goal of this paper is to give such a bijection which has useful
implications in terms of the geometry, topology and combinatorics of
affine Grassmannians and Bruhat order.  The family of partitions in
the image of this map is not the most obvious one: partitions whose
parts are all in the set of exponents.  Instead, we map $\minreps$ to
a family of colored partitions we call \textit{affine partitions}
using a canonical factorization into \textit{segments}.  The segments
are determined in general by the minimal length coset representatives
in a corresponding finite Weyl group, hence there are only a finite
number of them.  In the simplest cases, the segments are the $W$-orbit
of the special generator $s_{0}$ in $\minreps$ acting on the left.  In
other cases, the map works best if we use a smaller set of segments
and their images under an automorphism of the Dynkin diagram.

Using our bijection between $\minreps$ and affine partitions, there
are three natural partial orders on affine partitions. Bruhat order
and the left weak order on $\minreps$ are inherited from
$\widetilde{W}$.  Thus, the affine partitions also inherit these two
poset structures via the bijection.  In addition we will introduce a
generalization of Young's lattice on affine partitions which refines
the weak order and is refined by the Bruhat order.

In type $A$, Misra and Miwa \cite{Misra-Miwa} showed that $k$-core
partitions are in bijection with $\minreps$.  Erikson-Erikson
\cite{EE-98}and Lapoint-Morse \cite{Lapointe-Morse-2005} have shown
that the elements of $\minreps$ are in bijection not only with
$k$-core partitions, but also with $k$-bounded partitions, and skew
shapes with no long hooks.  Each of these partition bijections has
useful properties in terms of the geometry of affine Grassmannians.
Affine partition in type $A$ give a new perspective on these well
studied families.  From the point of view of affine partitions though,
type $A$ is harder than the other types
$B_{n},C_{n},D_{n},E_{6,7,8},G_{2},F_{4}$ because the weak order on
segments is the most complicated.  Therefore, type $A$ is covered last
though the reader interested only in type $A$ can skip past the other
type specific sections.

The segments and reduced factorizations have been used before in a
wide variety of other work
\cite{BME,BME-II,BME-III,EE-98,LLMS,LSS-2007,littig,Reiner-96} in
various types connecting affine partitions with lecture hall
partitions, Pieri type rules for the homology and cohomology of affine
Grassmannians, and hypergeometric identities.  However, none of the
previous work seems to address the complete set of affine Weyl groups
as we do in this article.  More details on previous work are given
after the definitions in Section~\ref{s:canonical}.

As an application of the theory of affine partitions and the
generalized Young's lattice, we will give a characterization of
rationally smooth Schubert varieties in the affine Grassmannians
$\loops$.  The smooth and rationally smooth Schubert varieties in
$\loops$ for all simple connected Lie Groups $G$ have recently been
characterized by the following theorem which was proved using
different techniques.  The description given here for the rationally
smooth Schubert varieties complements our original description and
relates these elements to the affine partitions.

\begin{theorem}\label{t:billey.mitchell}\cite{Billey-Mitchell}
Let $X_{w}$ be the Schubert variety in $\loops$ indexed by $w \in
\minreps$.  
\begin{enumerate}
\item $X_{w}$ is smooth if and only if $X_{w}$ is a closed parabolic
orbit.
\item $X_{w}$ is rationally smooth if and only if one of the following
conditions holds:
\begin{enumerate}
\item [a)] $X_w $ is a closed parabolic orbit.
\item [b)] The set $\{v \in \minreps: v\leq w \}$ is totally
ordered.
\item [c)]$W$ has type $A_{n}$ and $X_{w}$ is spiral (see
Section~\ref{s:type.a} for definition).
\item [d)]$W$ has type $B_3$ and
$w=s_{3}s_{2}s_{0}s_{3}s_{2}s_{1}s_{3}s_{2}s_{0}$.
\end{enumerate}
\end{enumerate}
\end{theorem}

By Theorem~\ref{t:billey.mitchell}, we will say $w$ is a \textit{cpo}
if $X_{w}$ is a closed parabolic orbit if and only if $X_{w}$ is
smooth.  In terms of Bruhat order, the cpo's can be identified as
follows.  Let $\widetilde{S}$ be a generating set for $\widetilde{W}$
and let $I_{w} = \{s \in \widetilde{S} : s \leq w \}$.  Then $w\in
\minreps$ is a cpo if and only if $sw \leq w$ for all $s \in I_{w}$.

For our characterization of rational smoothness, we will rely on the
following theorem due to Carrell and Peterson which requires a bit more
terminology.  It is known that the Poincar\'e polynomial of the
Schubert variety $X_{w}$ in $\loops$ is determined by
\[
P_{w}(t)= \sum t^{l(v)}
\]
where the sum is over all $v$ in $\minreps$ such that $v \leq w$ in
Bruhat order on $\widetilde{W}$.  See \cite{Kumar} for details.  We
say that a polynomial $F(t) = f_{0} + f_{1}t + f_{2}t^{2}+ \dots +
f_{k} t^{k} $ is \textit{palindromic} if $f_{i} = f_{k-i}$ for all
$0\leq i\leq k$.

\begin{theorem}\label{t:rationally.smooth}\cite{carrell94}
Let $X_{w}$ be the Schubert variety in $\loops$ indexed by $w \in
\minreps$.  Then $X_{w}$ is rationally smooth if and only if $P_{w}(t)$ is
palindromic.
\end{theorem}

In light of Theorem~\ref{t:rationally.smooth}, we will say $w\in
\minreps$ is \textit{palindromic} if and only if $P_{w}(t)$ is
palindromic if and only if $X_{w}$ is rationally smooth.  We will say
$w$ is a \textit{chain} if $\{v\in \minreps : v\leq w \}$ is a totally
ordered set.

The outline of the paper is as follows.  In Section~\ref{s:background}
we establish our basic notation and concepts we hope are familiar to
readers.  In Section~\ref{s:canonical}, the canonical factorization
into segments for elements in $\minreps$ is described for all Weyl
groups which motivates the definition of the affine partitions.  We
also state the main theorem giving the bijection from $\minreps$ to
affine partitions.  The type dependent part of the proof of the main
theorem is postponed until Sections~\ref{s:typeB}
through~\ref{s:type.a}.  In Section~\ref{s:palindromics}, we present a
new characterization of palindromic elements in terms of affine
partitions and generalized Young's lattice.

After posting the original version of this manuscript on arXiv.org, we
learned of the work of Andrew Pruett which also gives canonical
reduced expressions for elements in $\minreps$ for the simply laced
types and characterizes the palindromic elements \cite{Pruett}.

\section{Background}\label{s:background}

In this section we establish notation and terminology for Weyl groups,
affine Weyl groups and partitions.  There are several excellent
textbooks available which cover this material more thoroughly
including \cite{b-b,bou456,Hum} for (affine) Weyl groups and
\cite{Andrews-book,M1,ec1,ec2} for partitions.

Let $S=\{s_{1},\dots , s_{n} \}$ be the simple generators for $W$ and
let $s_{0}$ be the additional generator for $\widetilde{W}$.  Let $D$
be the Dynkin diagram for $\widetilde{W}$ as shown on
Page~\pageref{f:dynkin}.  Then the relations on the generators are
determined by $D$
\[
\left(s_{i} s_{j} \right)^{m_{ij}} =1
\]
where $m_{ij}=2$ if $i,j$ are not connected in $D$ and otherwise
$m_{ij}$ is the multiplicity of the bond between $i,j$ in $D$.  A
product of generators is \textit{reduced} if no shorter product
determines the same element in $\widetilde{W}$.

Let $\ell(w)$ denote the \textit{length} of $w \in \widetilde{W}$ or the
length of any reduced expression for $w$.  The Bruhat order on
$\widetilde{W}$ is defined by $v \leq w$ if given any reduced expression
$w=s_{a_{1}}s_{a_{2}}\cdots s_{a_{p}}$ there exists a subexpression
for $v$.  Therefore, the cover relation in Bruhat order is defined by 
\[
w \text{ covers } v \iff v = s_{a_{1}}s_{a_{2}}\cdots
\widehat{s_{a_{i}}} \cdots s_{a_{p}} \text{ and } \ell(w) = \ell(v) +1.
\]

A \textit{partition} is a weakly decreasing sequence of positive
integers of finite length.  By an abuse of terminology, we will also
consider a partition to be a weakly decreasing sequence of
non-negative integers with a finite number of positive terms.  A
partition $\lambda=(\lambda_{1}\geq \lambda_{2}\geq \lambda_{3} \geq
\ldots) \ $ is often depicted by a Ferrers diagram which is a left
justified set of squares with $\lambda_{1}$ squares on the top row,
$\lambda_{2}$ squares on the second row, etc.  For example,
\[
(7,5,5,2) \cong 
\tableau{{} &   {} &   {} &   {} &   {} &   {} &   {} \\
{} &   {} &   {} &   {} &   {} &   \\
{} &   {} &   {} &   {} &   {} &   \\
{} &   {} 
}
\]
The values $\lambda_{i}$ are called the \textit{parts} of the
partition.

Young's lattice on partitions is an important partial order determined
by containment of Ferrers diagrams \cite{M1}. In other words,
$\mu \subset \lambda $ if $\mu_{i} \leq \lambda_{i}$ for all $i\geq
0$.  For example, $(5,5,4) \subset (7,5,5,2)$ in Young's lattice.
Young's lattice is a ranked poset with rank function determined by the
\textit{size of the partition}, denoted
\[
|\lambda | = \sum \lambda_{i}.
\]

Young's lattice appears as the closure relation on Schubert varieties
in the classical Grassmannian varieties \cite{Fulton-book}. For the
isotropic Grassmannians of types $B,C,D$, the containment relation on
Schubert varieties is determined by the subposet of Young's lattice on
\textit{strict partitions}, i.e. partitions of the form
$(\lambda_{1}>\lambda_{2}>\dots > \lambda_{f})$.  This fact follows
easily from the signed permutation notation for the Weyl group of
types $B/C$.  Lascoux \cite{Lascoux-99} has shown that Young's lattice
restricted to $n+1$-core partitions characterizes the closure relation
on Schubert varieties for affine Grassmannians in type $A$.  See
Section~\ref{s:type.a} for more details.

\section{Canonical Factorizations and Affine Partitions}\label{s:canonical} 

In this section, we will identify a canonical reduced factorization
$r(w)$ for each minimal length coset representative $w \in \minreps$.
The factorizations will be in terms of \textit{segments} coming from
quotients of parabolic subgroups of $\widetilde{W}$.  We will use the
fact about Coxeter groups that for each $w \in \widetilde{W}$ and each
parabolic subgroup $W_{I}=\langle s_i \given i \in I \rangle$ there exists a unique factorization of $w$
such that $w=u \cdot v$, $\ell(w)= \ell(u) + \ell(v)$, $u$ is a
minimal length element in the coset $uW_{I}$ and $v \in W_{I}$
\cite[Prop. 2.4.4]{b-b}.  We will use the notation
$u=u_{I}(w)$ and $v=v_{I}(w)$ in this unique factorization of $w$.
Let $W^{I}$ denote the minimal length coset representatives for
$W/W_{I}$.

Consider the Coxeter graph of an affine Weyl group $\widetilde{W}$ as
labeled on Page~\pageref{f:dynkin}.  The special generator $s_{0}$ is
connected to either one or two elements among $s_{1},\dots, s_{n}$.
If $s_{0}$ is connected to $s_{1}$, call $\widetilde{W}$ a
\textit{Type I Coxeter group}; types $A,C,E,F,G$.  If $s_{0}$ is not
connected to $s_{1}$, then there is an involution on the Coxeter graph
for $\widetilde{W}$ interchanging $s_{0}$ and $s_{1}$ and fixing all
other generators.  Call these $\widetilde{W}$ \textit{Type II Coxeter
groups}; types $B,D$.

Let $\widetilde{W}$ be a Type I Coxeter group, then the parabolic
subgroup generated by $S=\{s_{1},s_{2}, \ldots , s_{n} \}$ is the
finite Weyl group $W$.  Let $J \subset S$ be the subset of generators
that commute with $s_{0}$; in particular $s_{1} \not \in J$ using our
labeling of the generators.  Then since $W$ is finite, there are a
finite number of minimal length coset representatives in $W^{J}$.  For
each $j\geq 0$, if there are $k$ elements in $W^{J}$ of length $j$,
label these \textit{fragments} by
\[
\fragment{1}(j),\ldots, \fragment{k}(j).
\]
Appending an $s_{0}$ onto the right of each fragment, we obtain
elements in $\minreps$ called \textit{segments}.  In particular, for
each fragment $\fragment{i}(j)$, fix a reduced expression
$\fragment{i}(j) = s_{a_{1}}s_{a_{2}}\cdots s_{a_{j}}$, and set
\[
\segment{i}{}(j+1) = s_{a_{1}}s_{a_{2}}\cdots s_{a_{j}} \cdot s_{0} \in \minreps.
\]
Now, assume $w \in \minreps$ and $w \neq \id$.  Let $w'=ws_{0}$.  Then
the unique reduced factorization $w'=u_{S}(w')\cdot v_{S}(w')$ has the
property that $u_{S}(w')$ is in $\minreps$ and $v_{S}(w') \in W$.  In
fact, $v_{S}(w') \in W^{J}$ since $w \in \minreps$ and all the
generators in $J$ commute with $s_{0}$.  Hence, $v_{S}(w') =
\fragment{i}{}(j)$ for some $i,j$ so multiplying on the right by
$s_{0}$ we have $w = u_{S}(w')\segment{i}{}(j+1)$.  By induction
$u_{S}(w')$ has a reduced factorization into a product of segments as
well.  Therefore, each $w \in \minreps$ has a canonical reduced
factorization into a product of segments, denoted
\begin{equation}\label{e:unique.expression.I}  
r(w) = \ \segment{i_{f}}{}(\lambda_{f}) \cdots
\segment{i_{3}}{}(\lambda_{3}) \segment{i_{2}}{}(\lambda_{2})
\segment{i_{1}}{}(\lambda_{1}),
\end{equation}
for some $f\geq 0$ and $l(w) = \sum_{j=1}^{f} \lambda_{j}$.  Note, $w
\in \minreps$ may have other reduced factorizations into a product of
segments, however, $r(w)$ is unique in the following sense.

\begin{lemma}\label{l:unique.factorization}
For $w\in \minreps$, the canonical reduced factorization $r(w)=
\segment{i_{f}}{}(\lambda_{f}) \cdots \\
\segment{i_{2}}{}(\lambda_{2})
\segment{i_{1}}{}(\lambda_{1})$ is the unique reduced factorization of
$w$ into a product of segments such that every initial product
\[
\segment{i_{f}}{}(\lambda_{f}) \cdots
\segment{i_{d+1}}{}(\lambda_{d+1}) \segment{i_{d}}{}(\lambda_{d}) 
 \hspace{.5in} 1\leq d\leq f
\]
is equal to $r(u)$ for some $u \in \minreps$.  Furthermore, every
consecutive partial product
\[
\segment{i_{d}}{}(\lambda_{d})\segment{i_{d-1}}{}(\lambda_{d-1})
\dotsb \segment{i_{c+1}}{}(\lambda_{c+1})
\segment{i_{c}}{}(\lambda_{c})  \hspace{.5in} 1\leq c\leq d\leq f
\]
is equal to $r(v)$ for some $v \in \minreps$ .
\end{lemma}

\begin{proof}
The first claim follows by induction from the uniqueness of the
factorization $w=u_{S}(w')\cdot v_{S}(w')s_{0}$.  To prove the second
claim, it is enough to assume $c=1$ by induction.  Let $v =
\segment{i_{d}}{}(\lambda_{d}) \dotsb \segment{i_{1}}{}(\lambda_{1})$
as an element of $\widetilde{W}$.  Then $v\in \minreps$ since $w \in
\minreps$, and the product $\segment{i_{d}}{}(\lambda_{d}) \dotsb
\segment{i_{2}}{}(\lambda_{2}) \segment{i_{1}}{}(\lambda_{1})$ must be
a reduced factorization of $v$ since $r(w)$ is a reduced
factorization.  Assume by induction on the number of segments in the
product that $\segment{i_{d}}{}(\lambda_{d}) \dotsb
\segment{i_{2}}{}(\lambda_{2}) = r(u)$ for some $u \in \minreps$.
Then, by the uniqueness of the factorization $r(v)$, we must have
$\segment{i_{d}}{}(\lambda_{d}) \dotsb \segment{i_{1}}{}(\lambda_{1})
=r(v)$.
\end{proof}

\begin{rem}
Observe that the canonical factorization into segments used above
extends to any Coxeter system $(W,S)$ with $s_{0}$ replaced by any
$s_{i} \in S$ and $J$ replaced by the set $\{s_{j} : j \neq i \text{
and } s_{i}s_{j}=s_{j}s_{i}\}$.  However, for types $B$ and $D$,
Lemma~\ref{l:pairs}(\ref{l:weak.order.and.segmetns}) and
Theorem~\ref{t:partition.bijection} don't hold using this
factorization.  By using the involution interchanging $s_{0}$ and
$s_{1}$ in Type II Weyl groups, we can identify another canonical
factorization with shorter segments and all the nice partition
properties as with Type I Weyl groups.
\end{rem}

Assume $\widetilde{W}$ is a Type II affine Weyl group with generators
labeled as on Page~\pageref{f:dynkin}.  Let $J=\{s_{2},s_{3},\dotsc ,
s_{n} \}$.  Note, the parabolic subgroups generated by
$S=\{s_{1},s_{2}, s_{3},\ldots , s_{n} \}$ and $S'=\{s_{0},s_{2},
s_{3},\ldots , s_{n} \}$ are isomorphic finite Weyl groups.  Since
$\widetilde{W}_{S'}$ is a finite Weyl group,
$(\widetilde{W}_{S'})^{J}$ is finite.  For each $j\geq 0$, if there
are $k$ elements in $(\widetilde{W}_{S'})^{J}$ of length $j$, label
these 0-\textit{segments} by
\begin{equation}\label{e:segments.II.0}
\segment{1}{0}(j),\ldots, \segment{k}{0}(j)
\end{equation}
and fix a reduced expression for each one.  
Similarly, there are a finite number of elements in $W^{J} =
(\widetilde{W}_{S})^{J}$.  For each $j\geq 0$, if there are $k$
minimal length elements of $W^{J}$ of length $j$, label these
1-\textit{segments} by
\begin{equation}\label{e:segments.II.1}
\segment{1}{1}(j),\ldots, \segment{k}{1}(j).
\end{equation}
We will assume each 0,1-pair of segments is labeled consistently
so $\segment{i}{1}(j)$ and $\segment{i}{0}(j)$ have reduced
expressions that differ only in the rightmost generator.  

By construction, every $\segment{i}{0}(j)$ is a minimum length coset
representative for $\quot = \widetilde{W}/W_{S}$ and every
$\segment{i}{1}(j)$ is a minimum length coset representative for
$\widetilde{W}/W_{S'}$.  Let $w \in \minreps$, then the unique reduced
factorization $w=u_{S'}(w)\cdot v_{S'}(w)$ has the property that
$u_{S'}(w) \in \widetilde{W}^{S'}$ and $v_{S'}(w) \in W_{S'}$.  In
fact, since $w \in \minreps$, then $v_{S'}(w) \in (W_{S'})^{J}$ so
$v_{S'}(w) = \segment{i}{0}(j)$ for some $i,j$.  Similarly, if $y \in
\widetilde{W}^{S'}$, then $y$ has a unique factorization $y = u_{S}(y)
\cdot v_{S}(y)$ where $u_{S}(y) \in \minreps$ and $v_{S}(y) =
\segment{i}{1}(j)$ for some $i,j$. Therefore, by induction each $w \in
\minreps$ has a canonical reduced factorization into a product of
alternating 0,1-segments, denoted
\begin{equation}\label{e:unique.expression.II}  
r(w) = \ \cdots \segment{c_{3}}{0}(\lambda_{3})
\segment{c_{2}}{1}(\lambda_{2}) \segment{c_{1}}{0}(\lambda_{1}),
\end{equation}
which is unique in the sense of Lemma~\ref{l:unique.factorization} but
where all the consecutive partial products correspond with minimum
length coset representatives in $\minreps$ or $\widetilde{W}^{S'}$
according to their rightmost factor.  Note, each $r(w)$ is the product
of a finite number of segments, say $f$ of them, and $l(w) =
\sum_{j=1}^{f} \lambda_{j}$.  Note further that the subscripts in
\eqref{e:unique.expression.II} are forced to start with 0 on the right
and then alternate between 0 and 1.  Hence the subscripts can easily
be recovered if we omit them and simply use the same notation as in
\eqref{e:unique.expression.I}.

In both Type I and Type II affine Weyl groups,
Lemma~\ref{l:pairs}(\ref{l:weak.order.and.segmetns}) below shows that
if $r(w)$ factors into segments of lengths $\lambda_{1},
\lambda_{2},\dotsc , \lambda_{f}$ as in \eqref{e:unique.expression.I}
or \eqref{e:unique.expression.II}, then the sequence of numbers
$(\lambda_{1},\lambda_{2},\dots)$ is a partition of $\ell(w)$ i.e.
$\lambda_{1}\geq \lambda_{2}\geq \lambda_{3}\geq \ldots \geq 0 $ and
$|\lambda | = \sum \lambda_{i} = \ell(w)$.  If the segments all have
unique lengths, then we can recover $w$ from the partition $\lambda$
by multiplying the corresponding segments in reverse order.  However,
when there are multiple segments of the same length we will need to
allow the parts of the partitions to be ``colored'' to be able to
recover $w$ from the colored partitions.  The colors of the parts of a
colored partition will be denoted by superscripts.  For example,
$(5^{1},5^{1},4^{2},3^{6},1^{2})$ corresponds with the partition
$(5,5,4,3,1)$ and the exponents determine the coloring of this
partition.  Colored partitions are only needed in types $A,D,E,F$.
Some colored partitions cannot occur in each type.  The rules for
determining the allowed colored partitions come from identifying which
products of pairs of segments are minimal length coset representatives
and which are not.

We will say that $(i^{a},j^{b})$ is an \textit{allowed pair} if
$\segment{a}{}(i) \cdot \segment{b}{}(j) \in \minreps$ and
$\ell(\segment{a}{}(i) \cdot \segment{b}{}(j)) = i+j$.  The following
two statements describe how segments and allowed pairs relate to the
left weak order on $\minreps$.

\begin{lemma} \label{l:pairs} If $(i^{a},j^{b})$ is an allowed pair,
then the following hold:
\begin{enumerate} 
\item \label{l:weak.order.and.segmetns} We have $\segment{a}{}(i) \leq
\segment{b}{}(j) \in \minreps$ in the left weak order on $\minreps$.
In particular, $i\leq j$ and $i=j$ implies $a=b$.
\\
\item \label{l:left.order.allowed} If $u \in \minreps$, $u \neq id$,
and $u < \segment{a}{}(i)$ in left weak order, then $u$ is a segment
itself, say $u=\segment{c}{}(h)$, and $(h^{c},j^{b})$ is also an allowed
pair.
\end{enumerate}
\end{lemma}

Note, $(h^{c},i^{a})$ may or may not be an allowed pair if
$\segment{c}{}(h) < \segment{a}{}(i)$.

\begin{proof}
The first statement in the lemma follows by observation in each type
once we have identified the segments in later sections.  In Type II
affine Weyl groups, the statement should be interpreted as: If
$(i^{a},j^{b})$ is an allowed pair, then $\segment{a}{0}(i) \leq
\segment{b}{0}(j)$ and $\segment{a}{1}(i) \leq \segment{b}{1}(j)$ in
left weak order.  To prove the first part of the second statement,
note that the set of segments form a lower order ideal in left weak
order.  Since $(i^{a},j^{b})$ is an allowed pair and $\segment{c}{}(h)
< \segment{a}{}(i)$ in left weak order, then $\segment{c}{}(h) \leq
\segment{b}{}(j)$ in left weak order also by (1).  Furthermore,
$\segment{c}{}(h)\cdot \segment{b}{}(j)$ is a right factor of
$\segment{a}{}(i) \cdot \segment{b}{}(j) \in \minreps$ so
$\segment{c}{}(h)\cdot \segment{b}{}(j)$ is reduced and in $\minreps$,
thus, $(h^{c},j^{b})$ is also an allowed pair.
\end{proof}

\begin{rem}
 It would be nice to have a type-independent proof of
Lemma~\ref{l:pairs}(\ref{l:weak.order.and.segmetns}) and to know to
what extent this statement holds for all Coxeter groups.   
\end{rem}

Let $\partitions$ be the set of \textit{affine colored partitions}
\textit{i.e.} the set of all partitions $ \lambda
=(\lambda_{1}^{c_{1}},\lambda_{2}^{c_{2}},\ldots,\lambda_{f}^{c_{f}})$
such that each consecutive pair
$(\lambda_{i+1}^{c_{i+1}},\lambda_{i}^{c_{i}})$ is an allowed pair for
$1\leq i<f$.  Note, the order of the consecutive pairs is backwards to
the order they appear in $\lambda$.  Observe that for every $w \in
\minreps$, the pairwise consecutive segments in $r(w)$ must correspond
with allowed pairs by Lemma~\ref{l:unique.factorization}.  Therefore,
we can define a map
\begin{equation}\label{e:partition.bijection}
\begin{array}{rcl}
\pi  : \minreps & \longrightarrow &\partitions\\
w  & \mapsto & \lambda 
\end{array}
\end{equation}
if $r(w)= \segment{c_{f}}{}(\lambda_{f}) \cdots
\segment{c_{2}}{}(\lambda_{2}) \segment{c_{1}}{}(\lambda_{1})$ and
$\lambda
=(\lambda_{1}^{c_{1}},\lambda_{2}^{c_{2}},\ldots,\lambda_{f}^{c_{f}})$.

\begin{theorem}\label{t:partition.bijection}
Let $W$ be any Weyl group and $\widetilde{W}$ be the corresponding
affine Weyl group.  Then $\pi : \minreps \rightarrow \partitions$ is a
length preserving bijection.
\end{theorem}

The theorem above is closely related to theorems of Lam, Lapoint,
Morse and Shimozono \cite{Lapointe-Morse-2005,LLMS} which have been
useful in the formulation of Pieri type rules for the cohomology ring
of the affine Grassmannian in type $A_{n}$.  The same bijection in
different language has been used by Eriksson-Eriksson \cite{EE-98} and
Reiner \cite{Reiner-96} to obtain partition identities related to
these bijections in types $B,C,D$.  See also Bousquet-M\'elou and
Eriksson \cite{BME,BME-II,BME-III} to relate $\minreps$ with the
\textit{lecture hall partitions}.  More recently
Lam-Shimozono-Schilling \cite{LSS-2007} and Littig \cite{littig} have
used similar factorizations in various types.
Our work provides a more general context in which the affine colored
partitions are related to $\minreps$ and the affine Grassmannians.

\begin{proof}
Observe that the map $\pi $ is automatically injective since
$r(v)=r(w)$ implies $v=w$.  Furthermore, the inverse map sending
\[
\lambda =(\lambda_{1}^{c_{1}} \geq \lambda_{2}^{c_{2}} \geq
\lambda_{3}^{c_{3}} \geq \dotsc ) \mapsto 
\begin{cases}
\ldots \segment{c_{3}}{}(\lambda_{3}) \segment{c_{2}}{}(\lambda_{2})
\segment{c_{1}}{}(\lambda_{1}) &	\text{Type I }\\
\ldots \segment{c_3}{0}(\lambda_{3})
\segment{c_2}{1}(\lambda_{2})\segment{c_1}{0}(\lambda_{1}) &
\text{Type II }
\end{cases}
\]
determines a well defined expression $\pi^{-1}(\lambda) $ in
$\widetilde{W}$.  If $\pi^{-1}(\lambda) \in \minreps$ and this
expression is reduced, then it must be $r(\pi^{-1}(\lambda))$ by
Lemma~\ref{l:unique.factorization}.  Therefore, after identifying the
segments and allowed pairs in each type, the theorem will follow if we
prove
\begin{quote}
\item For each $k\geq 0$, the number of partitions of $k$ in
$\partitions$ is equinumerous to the number of elements in $\minreps$
of length $k$.
\end{quote}
This statement can be proved via a partition identity equating Bott's
formula \eqref{e:bott.formula} and the rank generating function for
$\partitions$ for all types except type $A$.

This verification occurs in Theorem~\ref{t:surjective.Bn} for type
$B$, Theorem~\ref{t:surjective.Cn} for type $C$, and
Theorem~\ref{t:surjective.Dn} for type $D$.  For type $G_{2}$, this
partition identity is easy to check by hand.  For types
$E_{6},E_{7},E_{8},F_{4}$, computer verification of the identity can
be used as discussed in Section~\ref{s:exceptional}.

In type $A_{n}$ for $n\geq 2$, the generating function for allowed
partitions is not as easy to write down in one formula simultaneously
for all $n$.  Therefore, surjectivity is proved in
Theorem~\ref{t:surjective.An} using the $n+1$-core partitions.
\end{proof}

For emphasis, we state the following corollary of
Theorem~\ref{t:partition.bijection} which is a useful tool for the
applications.  As opposed to multiplication of generators, the
corollary says that reduced multiplication of segments is a ``local
condition''.

\begin{corollary}\label{c:products}
Any product of segments
$\segment{c_{1}}{}(j_{1})\segment{c_{2}}{}(j_{2})\dotsb
\segment{c_{k}}{}(j_{k})$ is equal to $r(v)$ for some $v \in \minreps$
if and only if for each $1\leq i<k$ the pair
$(j_{i}^{c_{i}},j_{i+1}^{c_{i+1}})$ is an allowed pair.
\end{corollary}

Given an affine partition in $\partitions$, say a corner is
\textit{$\partitions$-removable} if the partition obtained by removing this corner
in the Ferrers diagram leaves a partition that is still in
$\partitions$.  The set of $\partitions$-removable corners for any partition will
depend on the affine Weyl group type. 

In types $B,C,G_{2}$, we will show that the segments have unique
lengths so $\partitions$ is a subset of all partitions with no colors
necessary.  It is interesting to note the relationship between Bruhat
order on $\minreps$ and the induced order from Young's lattice on
$\partitions$.  The corollary below shows that Bruhat order on
$\minreps$ contains the covering relations in Young's lattice
determined by $\partitions$-removable corners.

\begin{corollary}\label{c:bruhat.shape.covers}
Let $\partitions$ be  the  set  of     affine partitions  in  types
$B_{n},C_{n}$  or  $G_{2}$.  If   $\lambda,  \mu \in  \partitions$ and
$\lambda$ covers $\mu$ in  Young's lattice, then $\pi^{-1} (\lambda )$
covers $\pi^{-1} (\mu )$ in the Bruhat order on $\minreps$.
\end{corollary}

\begin{proof}
In types $B_{n}$, $C_{n}$, and $G_{2}$, the segments form a chain in
the left weak order, see Equations~\eqref{e:Bsegments.1},
\eqref{e:Csegments.1}, and \eqref{e:segments.g2}.  If $\lambda, \mu
\in \partitions$ and $\lambda$ covers $\mu$ in Young's lattice, then
$\mu$ is obtained by deleting one outside corner square from
$\lambda$.  Thus $\segment{}{}(\mu_{g})\cdots
\segment{}{}(\mu_{2})\segment{}{}(\mu_{1})$ is obtained from
$\segment{}{}(\lambda_{f})\cdots
\segment{}{}(\lambda_{2})\segment{}{}(\lambda_{1})$ by striking out
one generator at the beginning of a segment.  This expression will be
a reduced expression for a minimal length coset representative
precisely when the corresponding partition satisfies the conditions to
be in $\partitions$ by Corollary~\ref{c:products}.  Given that $\mu
\in \partitions$, then $\pi^{-1} (\lambda ) =
\segment{}{}(\lambda_{f})\cdots
\segment{}{}(\lambda_{2})\segment{}{}(\lambda_{1}) >
\segment{}{}(\mu_{g})\cdots \segment{}{}(\mu_{2})\segment{}{}(\mu_{1})
= \pi^{-1} (\mu )$ and $\ell(\pi^{-1} (\lambda )) = \ell(\pi^{-1}(\mu
)) +1$.  So, $\pi^{-1}(\lambda )$ covers $\pi^{-1}(\mu )$ in the
Bruhat order on $\minreps$.
\end{proof}

For types $A$,$D$,$E$, and $F$, we define a generalization of Young's
lattice on colored partitions as follows.  First, a colored part
$j^{c}$ covers another part $(j-1)^{d}$ if $\segment{c}{}(j)$ covers
$\segment{d}{}(j-1)$ in left weak order on $\widetilde{W}$.  Second, a
colored partition $\lambda =
(\lambda_{1}^{c_{1}},\lambda_{2}^{c_{2}},\ldots) \in \partitions$
covers $\mu = (\mu_{1}^{d_{1}},\mu_{2}^{d_{2}},\ldots) \in
\partitions$ if $\lambda$ and $\mu$ agree in all but one part indexed
by $j$, and $\lambda_{j}^{c_{j}} $ covers $\mu_{j}^{d_{j}}$ in the
partial order on colored parts.

\begin{corollary}\label{c:gen.young.covers}
If $\lambda, \mu \in \partitions$ and $\lambda$ covers $\mu$ in the
generalized Young's lattice, then $\pi^{-1} (\lambda )$ covers
$\pi^{-1} (\mu )$ in the Bruhat order on $\minreps$.
\end{corollary}

\begin{rem}\label{r:converse}
The converse to Corollary~\ref{c:gen.young.covers} does not hold.  See
Example~\ref{ex:bruhat.B}.  
\end{rem}

\begin{question}\label{open:k-core.bn}
Is there an alternative partition $\xi (w)$ to associate with each $w
\in \minreps$ so that $v<w$ in Bruhat order on $\minreps$ if and only
if $\xi (v) \subset \xi (w)$ outside of type $A$?  Recall, the core
partitions play this role in type $A$ by a a theorem of Lascoux
\cite{Lascoux-99}. 
\end{question}

\section{Palindromic Elements}\label{s:palindromics}

As a consequence of the bijection between $\minreps$ and affine
partitions from Theorem~\ref{t:partition.bijection}, the generalized
Young's lattice and Corollary~\ref{c:gen.young.covers}, we can observe
enough relations in Bruhat order to identify all palindromic elements
of $\minreps$ in terms of affine partitions.  For example, any affine
partition with two or more $\partitions$-removable corners cannot correspond with a
palindromic element since Bott's formula starts $1+t+\dotsb$ in all
types.  We recall, the palindromic elements have been recently
characterized in \cite{Billey-Mitchell} via the coroot lattice
elements.  In type $A_{n}$, there are two infinite families of
palindromics, first studied by the second author in
\cite{Mitchell-86}.  This alternative approach has given us additional
insight into the combinatorial structure of $\minreps$.

\begin{theorem}\label{t:palindromics} Assume $W$ is not of type
$A_{n}$ for $n\geq 2$ or $B_{3}$.  Let $w \in \minreps$ and say
$\pi(w) = \lambda$.  Then $w$ is palindromic if and only if the
interval $[id,w]$ in Bruhat order on $\minreps$ is isomorphic to the
interval $[\emptyset, \lambda ]$ in the generalized Young's lattice
and the interval $[\emptyset, \lambda ]$ is rank symmetric.
\end{theorem}

\begin{rem}\label{r:yb-nice}
Say $w\in \minreps$ is \textit{$YB$-nice} if the interval $[id,w]$ in
Bruhat order on $\minreps$ is isomorphic to the interval $[\emptyset,
\lambda ]$ in the generalized Young's lattice.  Say $w$ is
\textit{$YB$-palindromic} if $w$ is $YB-nice$ and the interval
$[\emptyset, \lambda ]$ is rank symmetric.  So outside of type $B_{3}$
and $A_{n}$ for $n\geq 2$, $YB-palindromic$ and palindromic are
equivalent.  In type $B_{3}$, there is one palindromic element which
is not $YB$-nice namely
$w=s_{3}s_{2}s_{0}s_{3}s_{2}s_{1}s_{3}s_{2}s_{0}$.  In type $A_{n}$
for $n\geq 2$ the spiral elements which are not closed parabolic
orbits are palindromic but not $YB$-nice.  See Section~\ref{s:type.a}.
\end{rem}

Typically the palindromic elements are indexed by one row shapes and
staircase shapes.  These elements correspond with chains and closed
parabolic orbits.  The proof of Theorem~\ref{t:palindromics} for the
infinite families will be stated and proved more explicitly in
Theorems~\ref{t:rationally.smooth.B},~\ref{t:rationally.smooth.C}
and~\ref{t:rationally.smooth.D}. After stating some general tools used
for the palindromy proofs, we prove Theorem~\ref{t:palindromics} for
the exceptional types below to demonstrate the technique of using
affine partitions.

Define the \textit{branching number} $b_{W}$ to be the smallest rank
in the Bruhat order on $\minreps$ with more than 1 element.  By Bott's
formula, $b_{W}\geq 2$.  It will be shown that if $w \in \minreps$ is
not palindromic, then symmetry of the Poincar\'e polynomial $P_{w}(t)$
always fails in the first $b_{W}-1$ coefficients for the exceptional
types.  For some of the other types, we must look further down in
Bruhat order.  

Let $m_{W}$ be the maximum number of coefficients one must check for
all $w\in \minreps$ to insure that $P_{w}(t)$ is not palindromic if
$w$ is not palindromic.  Hence, $m_{W}$ is defined to be the minimum
number so that $w \in \minreps$ is palindromic if and only if
$P_{w}(t)= \sum_{i=0}^{\ell(w)} a_{i}t^{i}$ and $a_{i}=a_{\ell(w)}$
for all $1\leq i\leq m_{W}$.  We will show that $m_{W}$ is always
bounded by $n$ but can be significantly smaller.

\begin{theorem}\label{t:mw}
The number $m_{W}$ for $\minreps$ is determined as follows:

\begin{equation}\label{e:mw}
\begin{array}{|l|c|}
\hline 
\text{type } & m_{W} \\
\hline \hline 
A_{1}, n=1 &	0  \\
A_{n}, n\geq 2 &	2 \\
B_{n}, n=3 &	2 \\   
B_{n}, n\geq 4 &	4 \\
C_{n}, n\geq 2 &	2 \\   
D_{n}, n\geq 4 &	n-2 \\ 
E_{6} &	3\\
E_{7} &	 4 \\
E_{8} &	6\\
F_{4} & 4\\
G_{2} & 4	\\
\hline 
\end{array}
\end{equation}

\end{theorem}

This theorem will be proved after introducing the segments in each
type.  

\begin{rem}
Computationally, the fact that $m_{W}$ is constant in most cases is
very useful. This means that an efficient algorithm exists to verify
that an element is not palindromic which does not require one to build
up the entire Bruhat interval below an element $w \in \minreps$ (which
takes an exponential amount of time in terms of the length of the
element).  For example, in type $F_{4}$ it suffices to choose a single
reduced expression for $w$ and consider subsequences with at most 4
generators removed.  This leads to an $O(\ell(w)^{4})$ algorithm.
\end{rem}

\begin{rem}\label{}
The number $m_{W}$ bounds the degree of the first nontrivial
coefficient in the Kazhdan-Lusztig polynomial $P_{id,w}(t)$ for $w \in
\minreps$ by a theorem of Bjorner and Ekedahl \cite{bjorner-ekedahl}.
\end{rem}

We identify the elements of $\minreps$ with their corresponding affine
partition by Theorem~\ref{t:partition.bijection}.  Therefore, the
Bruhat order and the left weak order extends to affine partitions in
addition to the generalized Young's lattice.  Furthermore, we will
abuse notation and denote a colored partition
$(\lambda_{1}^{c_{1}},\dotsc, \lambda_{k}^{c_{k}})$ simply by
$\lambda$ or $(\lambda_{1},\dotsc, \lambda_{k})$ when the particular colors are not
essential to the argument.  

Given an affine partition $\lambda$, we will say $\lambda$ is a
\textit{thin partition} if the interval $[\emptyset, \lambda]$ is rank
symmetric in the first and last $b_{W}$ ranks in the generalized
Young's lattice.  Note, thin is a necessary condition for palindromy
by Corollary~\ref{c:gen.young.covers}.  Furthermore, the following
lemma shows that any affine partition whose smallest $k$ parts forms a
partition that is not thin cannot itself be thin.

\begin{lemma}\label{l:thin.factors}
If a colored partition $(\lambda_{1}^{c_{1}},\dotsc,
\lambda_{k}^{c_{k}}) \in \partitions$ is thin then
$(\lambda_{2}^{c_{2}},\dotsc, \lambda_{k}^{c_{k}}) \in \partitions$ is
thin.
\end{lemma}

\begin{proof}
Suppose $\lambda'=(\lambda_{2}^{c_{2}},\dotsc, \lambda_{k}^{c_{k}})$
is not thin.  Then there exist at least two affine colored partitions
$\mu,\nu$ below $\lambda '$ in generalized Young's lattice such that
$|\mu | =|\nu | > |\lambda '|-b_{W}$.  By definition of generalized
Young's lattice we must have $\mu_{1},\nu_{1} \leq
\lambda_{2}^{c_{2}}$ in left weak order.  Therefore, by
Lemma~\ref{l:pairs}, both $(\mu_{1},\lambda_{1}^{c_{1}})$ and
$(\nu_{1},\lambda_{1}^{c_{1}})$ must be allowed pairs.  Therefore,
concatenating $\lambda_{1}^{c_{1}}$ on the front of $\mu,\nu $ gives
two partitions $\mu ' = (\lambda_{1}^{c_{1}},\mu)$ and $\nu ' =
(\lambda_{1}^{c_{1}},\nu)$ in $\partitions$ such that $\mu ' , \nu '
\leq \lambda $ in generalized Young's lattice and $|\mu' | =|\nu' | >
|\lambda|-b_{W}$, contradicting the fact that $\lambda$ is thin.
Hence, $\lambda '$ must also be thin.
\end{proof}

Observe that $\lambda$ may be thin and still cover two or more
elements in Bruhat order; exactly one of the covering relations occurs
in the generalized Young's lattice in this case.  (See
Example~\ref{ex:bruhat.B}.)  Using similar logic to the proof of
Lemma~\ref{l:thin.factors}, we get the following statement.  In Bruhat
order, if a colored partition $\lambda = (\lambda_{1},\dotsc,
\lambda_{k})$ (suppressing the colors) covers two partitions $\mu,\nu
$ both with largest part at most $\lambda_1$ in left weak order, then
any affine partition of the form $ (\gamma_{1},\dotsc,
\gamma_{j},\lambda_{1},\dotsc, \lambda_{k}) = \gamma .\lambda $ also
covers two partitions $\gamma .\mu$ and $\gamma .\nu$.  Therefore,
define an \textit{extra thin partition} $\lambda$ to be an affine
partition that is thin and such that there exists at most one affine
partition $\mu$ covered by $\lambda$ in Bruhat order such that the
largest part of $\mu$ is less than or equal to the largest part of
$\lambda$ in left weak order.  Note, extra thin is a necessary
condition for palindromy.  In fact, this gives an analog of
Lemma~\ref{l:thin.factors} for extra thin partitions.

\begin{lemma}\label{l:extra.thin.factors}
If a colored partition $(\lambda_{1}^{c_{1}},\dotsc,
\lambda_{k}^{c_{k}}) \in \partitions$ is extra thin then the colored
partition
$(\lambda_{2}^{c_{2}},\dotsc, \lambda_{k}^{c_{k}})
 \in
\partitions$ is extra thin.
\end{lemma}

Now we restrict our attention to the exceptional types.  Determining
the palindromic elements in $\minreps$ for the exceptional types
follows easily from the affine partitions, the generalized Young's
lattice and Lemma~\ref{l:exceptional.finite}.  For each type there are
relatively few palindromics: 8 in $G_{2}$, 9 in $F_{4}$, and 11 in
$E_{6,7,8}$.  In types $G_{2}$ and $F_{4}$ the palindromics are
precisely the \textit{chains}: elements whose principal lower order
ideal is a chain in Bruhat order so $P_{w}(t)=1+t+\dots +
t^{\ell(w)}$.  In each of types $E_{6,7,8}$, in addition to the chains
there is one palindromic whose Poincar\'e polynomial is similar to the
longest single row in type $D$.  Plus in types $E_{6,7}$ there are 2
additional Poincar\'e polynomials that occur.  In $E_{6}$, two
elements have the Poincar\'e polynomial
\[
({t}^{10}+{t}^{9}+{t}^{8}+{t}^{6}+{t}^{5}+{t}^{4}+{t}^{2}+t+1) (
t^{6}+t^{3}+1).
\]
In $E_{7}$, one element has the Poincar\'e polynomial
\[
({t}^{21}+{t}^{16}+{t}^{14}+{t}^{12}+{t}^{9}+{t}^{7}+{t}^{5}+1) (
t^{6}+ t^{5}+ t^{4}+ t^{3}+ t^{2}+ t +1)
\]

\begin{lemma}\label{l:exceptional.finite}
For all exceptional types, every affine partition with 7 or more parts
is not palindromic.  In fact,  $m_{W}\leq (b_{W}-1)$.
\end{lemma}

\begin{proof}

Recall that by Lemma~\ref{l:pairs}, every non-empty affine partition
covers at least one element in the generalized Young's lattice by
removing a corner of the smallest part.  So by
Corollary~\ref{c:gen.young.covers}, given an affine partition $\mu $
with 7 or more parts, we only need to show there exists one element
$\nu$ below $\mu$ in Bruhat order which is not comparable to $\mu$ in
generalized Young's lattice and $|\mu|-|\nu| < b_{W}$.

By definition, any affine partition which is not extra thin covers two
or more elements in Bruhat order so we can restrict our attention to
the extra thin elements.  The following observations can be made in
the exceptional types with computer assistance:
\begin{enumerate}
\item The only extra thin affine partitions with 7 parts in any
exceptional type occur in $E_6, E_7$ and $G_2.$

\item In $E_6, E_7$ and $G_2$, every extra thin affine partition with
7 parts has largest 4 parts $(j^{c},j^{c},j^{c},j^{c})$ (repeats at
least 4 times).  Furthermore, if $k^{d}$ is any part such that
$(j^{c},k^{d})$ is any allowed pair, then either $j^{c} = k^{d}$ or
the affine partition $(k^{d},j^{c})$ has two $\partitions$-removable
corners.

\item The repeated parts $j^{c}$ occurring in every extra thin affine
partition with 7 parts described in Part (2) have the property that
below $(j^{c},j^{c},j^{c},j^{c})$ in Bruhat order there exists some
affine partition $\lambda$ whose parts are all weakly larger than
$j^{c}$ in the left weak order and $|\lambda|>4j -b_{W}$.
\end{enumerate}

Using these observations we complete the proof.  In types $E_{8}$ and
$F_{4}$ every affine partition with 7 or more parts is not extra thin
by the first observation and Lemma~\ref{l:extra.thin.factors}.  So,
assume the type is $E_6, E_7$ or $G_2$.  By the second observation,
the only colored part which is allowed to extend an extra thin affine
partition with 7 or more parts to another extra thin affine partition
is another copy of its largest part.  Therefore, every extra thin
affine partition $\mu$ with at least 7 parts has its largest part
repeated at least 4 times.  Say $\mu = (j^{c},j^{c},j^{c},j^{c},
\mu_{5}^{c_{5}},\dotsc, \mu_{k}^{c_{k}})$ and say $\lambda =
(\lambda_{1}^{d_{1}},\dotsc, \lambda_{p}^{d_{p}})$ is the affine
partition below $(j^{c},j^{c},j^{c},j^{c})$ in Bruhat order appearing
in the third observation.  Then $\nu = (\lambda_{1}^{d_{1}},\dotsc,
\lambda_{p}^{d_{p}}, \mu_{5}^{c_{5}},\dotsc, \mu_{k}^{c_{k}}) $ must
be an affine partition by Lemma~\ref{l:pairs} and the definition of
$\partitions$.  Furthermore, $\nu$ is not in the interval below $\mu$
in the generalized Young's lattice but $\nu <\mu $ in Bruhat order and
$|\nu| > |\mu|- b_{W}$.  Hence $\mu$ fails to be palindromic by the
$b_{W}-1$st coefficient.
\end{proof}

\begin{rem}
Looking closely at the data, one sees that in fact every affine
partition in $E_{7},E_{8}$ and $F_{4}$ with 7 or more parts covers at
least two elements in Bruhat order.  In $E_{6}$ and $G_{2}$, there
exist an infinite number of extra thin non-palindromic elements which
cover a single element in Bruhat order and only fail to be palindromic
on the second and fourth coefficient respectively.  In type $E_{6}$
there are two segments of length 12 which can repeat an arbitrary
number of times to get such elements.  In type $G_{2}$, there is a
unique segment of length 5, and the affine partition $(5,5,\dotsc,5)$
with $k\geq 2$ parts has this property.
\end{rem}

\begin{theorem}\label{t:exceptional.pal} In each exceptional type we
have:
\begin{enumerate}
\item Every affine partition $\lambda \in \partitions$ with 4 or more
parts corresponds with a non-palindromic $w \in \minreps$.  In fact,
the unique palindromic with 3 parts occurs in $E_{7}$.
\item A finite computer search over all extra thin affine partitions
with at most 3 parts suffices to identify all palindromic elements in
$\partitions$, or equivalently in $\minreps$.
\item We have $m_{W}=b_{W}-1$.
\end{enumerate}
\end{theorem}

\begin{proof}
Note, every element which is not extra thin fails to be palindromic by
depth $b_{W}-1$ by definition.  Therefore, by
Lemma~\ref{l:exceptional.finite}, we only needed to verify these
statements for extra thin affine partitions with at most 6 parts which
is efficient to check by computer.  In each case, it was verified that
every non-palindromic fails to be palindromic by depth $b_{W}-1$, and
furthermore, there is at least one non-palindromic element that is
palindromic up to the $(b_{W}-2)^{nd}$ coefficient.
\end{proof}

\begin{rem}
We note that the proof of Theorem~\ref{t:palindromics} in the
exceptional types now follows from Theorem~\ref{t:exceptional.pal}
simply by verifying that for every palindromic element the
corresponding Bruhat interval is isomorphic to the interval in
generalized Young's lattice.
\end{rem}



\section{Type $B$}\label{s:typeB} 

In this section we prove Theorem~\ref{t:partition.bijection} and
Theorem~\ref{t:palindromics} for type $B_{n}$, $n\geq 3$.  We begin by
identifying a family of partitions $\partitions(B_{n})$ with a length
preserving bijection to $\minreps$.  Then, we identify the segments in
type $B$ and the allowed pairs corresponding to these segments.  It
follows immediately from the list of allowed pairs that
$\partitions(B_{n})$ are the affine partitions.  Finally we use the
explicit description of the segments in this type to identify the
affine partitions corresponding with palindromic elements in
$\minreps$.

Let $\partitions(B_{n})$ be the set of partitions whose parts are
bounded by $2n-1$ and all the parts of length strictly less than $n$
are strictly decreasing.  Note, all parts in these partitions have the
same color.  Therefore, the generating function for such partitions is
\[
G_{\Bn}(x)=\frac{(1+x)(1+x^{2})\cdots (1+x^{n-1})}{(1-x^{n}) (1-x^{n+1}) \cdots (1-x^{2n-1})}
\]

\begin{lemma}\label{l:Bpartitions}
We have the following generating function identity with Bott's formula
from \eqref{e:bott.formula}:
\[
G_{\Bn}(x) = \frac{1}{(1-x)(1-x^{3})(1-x^{5})\cdots (1-x^{2n-1})}.
\]
\end{lemma}

\begin{proof}
Apply induction on $n$.
\end{proof}

Consider the Coxeter graph of $\widetilde{B}_{n}$ on
Page~\pageref{f:dynkin}.  In the language of
Section~\ref{s:canonical}, $\widetilde{B}_{n}$ is a Type II Coxeter graph.
Following \eqref{e:segments.II.0} and \eqref{e:segments.II.1}, for
$1\leq j \leq 2n-1$, set
\begin{equation}\label{e:Bsegments.1}
\segment{}{1}(j) = \begin{cases}
s_{j} \dots  s_{3}s_{2} s_{1} &	 1\leq j\leq n\\
s_{2n-j} \dotsb s_{n-1}s_{n} s_{n-1}\dotsb s_{4}s_{3}s_{2}  s_{1} &
n<j\leq 2n-1
\end{cases}
\end{equation}
Similarly, replacing all the $s_{1}$'s in $\segment{}{1}(j)$ with $s_{0}$'s, set 
\begin{equation}\label{e:Bsegments.2}
\segment{}{0}(j) = \begin{cases}
s_{0} &	 j=1\\
s_{j} \dots  s_{3}s_{2} s_{0} &	 1<j\leq n\\
s_{2n-j}  \dots  s_{n-1} s_{n} s_{n-1} \dots  s_{3}s_{2} s_{0} &
n<j\leq 2n-2\\
s_{0} s_{2} s_{3}\dots s_{n-1}s_{n} s_{n-1} \dots  s_{3}s_{2} s_{0} &	 j=2n-1.
\end{cases}
\end{equation}
Using the signed permutation description of the finite Weyl group of
type $B_{n}$, one can verify that these segments are all the minimal
length coset representatives for $W/W_{J}$ and
$\widetilde{W}_{S'}/W_{J}$ respectively.  See \cite[Chapter
8]{b-b} for a detailed description of this notation.

\begin{lemma}\label{l:commutation.B}
We have the following commutation rules for $s_{i} \cdot \segment{}{1}(j)$ for all
$1 \leq i \leq n $ and $1\leq j\leq 2n-1$:
\begin{align}\label{e:Bsegments}
s_{i} \cdot \segment{}{1}(j) = \begin{cases}
\segment{}{1}(j) \cdot s_{i} &	 1\leq j < i-1  \text{ or  }   2n-i<j \leq 2n-1\\
\segment{}{1}(j+1) &	 j = i-1 \text{ or } j = 2n -i-1\\
\segment{}{1}(j-1) &   j = i \text{ or  } j=2n-i\\
\segment{}{1}(j)\cdot s_{i+1} & 1\leq i < j < 2n-i-1,
\end{cases}
\end{align}
and for $i=0$ and $j=1$ we have 
\[
s_{0} \segment{}{1}(1) = s_{0}s_{1}=s_{1}s_{0}=\segment{}{1}(1) s_{0}.
\]
Similar commutation rules for $s_{i} \segment{}{0}(j)$ for all $0 \leq i \leq n $
and $1\leq j\leq 2n-1$ are obtained from \eqref{e:Bsegments} by
interchanging the roles of $s_{1}$ and $s_{0}$.  
\end{lemma}

\begin{proof}
These follow directly from the commutation relations among the
generators determined by the Dynkin diagram.
\end{proof}

\begin{lemma}\label{l:C.commutations.B}
We have the following product rules for segments in type $B_{n}$:
\begin{enumerate}
\item  For $1\leq j<n$, 
\[
\segment{}{1}(j) \segment{}{0}(j) =
\segment{}{1}(j-1)\segment{}{0}(j)\cdot s_{1}.
\]
\item  For $n\leq j \leq 2n-1$,
\[
\segment{}{1}(j+1) \segment{}{0}(j) =
\segment{}{1}(j)\segment{}{0}(j)\cdot s_{1}.
\]
\end{enumerate}
\end{lemma}

\begin{proof}
If $j=1$, then Statement 1 holds since $s_{0}$ and $s_{1}$ commute.
Assume $2\leq j<n$.  By Lemma~\ref{l:commutation.B}, we have
$\segment{}{1}(j) s_{i} = s_{i-1} \segment{}{1}(j) $ for all $2\leq i
\leq j$ so
\begin{align}
\segment{}{1}(j) \cdot \segment{}{0}(j) = & \segment{}{1}(j) \cdot s_{j}s_{j-1}\dotsb s_{2}s_{0}
\\
= & s_{j-1}\dotsb s_{2} s_{1} \cdot \segment{}{1}(j) \cdot s_{0} 
\\
= & \segment{}{1}(j-1) \cdot s_{j}s_{j-1}\dotsb s_{2}s_{1} \cdot s_{0}
\\
= &  \segment{}{1}(j-1) \cdot s_{j}s_{j-1}\dotsb s_{2}\cdot s_{0}  s_{1}
\\
= & \segment{}{1}(j-1) \cdot \segment{}{0}(j) \cdot s_{1}.
\end{align}


Statement 2 holds by a similar argument.  
\end{proof}

\begin{corollary}\label{c:C.all.commutations.B}
We have the following product rules for segments corresponding with
all pairs $(j,k) \not \in \partitions(B_{n})$ with $1\leq j,k\leq
2n-1$:
\[
\segment{}{1}(k) \cdot \segment{}{0}(j) = 
\begin{cases}
\segment{}{1}(j-1) \cdot \segment{}{0}(k) \cdot s_{1}  &	j \leq k < 2n-j\\
\segment{}{1}(j) \cdot \segment{}{0}(k-1) \cdot s_{1}  &	n \leq j < k \text{ or  }  j< 2n-j\leq k.
\end{cases}
\]
\end{corollary}

\begin{proof}
This follows from Lemma~\ref{l:commutation.B} and
Lemma~\ref{l:C.commutations.B}.
\end{proof}

Recall, $\partitions$ is the set of allowed partitions for $\minreps$
and $\pi : \minreps \longrightarrow \partitions$ was defined in
\eqref{e:partition.bijection}. 

\begin{theorem}\label{t:surjective.Bn}
In type $B_{n}$, we have $\partitions = \partitions(\Bn)$ and $\pi :
\minreps \longrightarrow \partitions$ is a length preserving
bijection.
\end{theorem}

\begin{proof}

We have already shown in the proof of
Theorem~\ref{t:partition.bijection} that $\pi$ is a length preserving
injection.  Since there is a length preserving bijection from
$\minreps$ to $\partitions(\Bn)$ by Lemma~\ref{l:Bpartitions}, we only
need to show that $\partitions \subset \partitions(\Bn)$ to prove the
theorem.  Therefore, we only need to show that the segments have
unique lengths between 1 and $2n-1$ and that all allowed pairs in type
$\Bn$ are strictly increasing if the larger part has length less than
$n$.  This follows directly from the definition of the segments and
Corollary~\ref{c:C.all.commutations.B}.
\end{proof}

The counting argument in the proof above determines the complete set
of allowed pairs.

\begin{corollary}\label{c:allowed.pairs.B}
The product $\segment{}{1}(j) \cdot \segment{}{0}(k) \in \minreps$ and
$\ell(\segment{}{1}(j) \cdot \segment{}{0}(k))=j+k$ if and only if
$1\leq j<k\leq 2n-1 $ or $n\leq j=k \leq 2n-1$.  
\end{corollary}

\begin{ex}\label{ex:labeling.B}
If $n=4$, the partition $(7,5,5,3,1) \in \partitions(B_{n})$
corresponds with $s_0\cdot s_{3}s_{2}s_{1} \cdot
s_{3}s_{4}s_{3}s_{2}s_{0} \cdot s_{3}s_{4}s_{3}s_{2}s_{1} \cdot
s_{0}s_{2}s_{3}s_{4}s_{3}s_{2}s_{0}$.  Pictorially, we have
\[
\tableau{
0 & 2 & 3 & 4 & 3 & 2 & 0\\
1 & 2 & 3 & 4 & 3 &  & \\
0 & 2 & 3 & 4 & 3 &  & \\
1 & 2 & 3 &  &  &  & \\
0 &  &  &  &  &  & }
\]
\end{ex}
\noindent where the corresponding reduced expression is read along the
rows from right to left, bottom to top.  Note, the 0's and 1's
alternate in the first column starting with a 0 in the top row since
$B_{n}$ is a Type II Coxeter group.  The other columns contain their
column number up to column $n$, for $n< j < 2n-1$ column $j$ contains
$2n-j$, and column $2n-1$ again alternates from 0 to 1 starting from
the top.

Recall from Corollary~\ref{c:bruhat.shape.covers}, that if both
$\lambda,\mu \in \partitions(B_{n})$ and $\lambda$ covers $\mu$ in
Young's lattice then $\pi^{-1}(\lambda)$ covers $\pi^{-1}(\mu)$ in
Bruhat order on $\minreps$.  However, below we give an example where the
converse does not hold.

\begin{ex}\label{ex:bruhat.B}
If $n=4$, the partition $(5,2,1) \in \partitions(B_{n})$ corresponds
with $s_{0} \cdot s_{2}s_{1} \cdot s_{3}s_{4}s_{3}s_{2}s_{0}$ which
covers $s_{0} s_{2} s_{3}s_{4}s_{3}s_{2}s_{0} = \segment{}{0}(7)$ by
deleting $s_{1}$.  Pictorially,
\[ (5,2,1) \cong
\tableau{
0 & 2 & 3 & 4 & 3 &  & \\
1 & 2 &  &  &  &  & \\
0 &  &  &  &  &  & }
\hspace{.2in}\mathrm{covers} \hspace{.2in}
(7) \cong
\tableau{
0 & 2 & 3 & 4 & 3 & 2 & 0\\
 &  &  &  &  &  & \\
 &  &  &  &  &  & }.
\]
Therefore, in Bruhat order $(5,2,1)$ covers $(7)$ even though
$(5,2,1)$ is not comparable to $(7)$ in Young's lattice.
\end{ex}

With more work, we can obtain all of the elements below $w \in
\minreps$ in Bruhat order by looking at $r(w)$, knocking out one
generator at a time, and using the commutation relations in
Lemma~\ref{l:commutation.B} and Corollary~\ref{c:C.all.commutations.B}
to identify a canonical segmented expression for the new element.  In
fact, in practice it is easy to see that knocking out most elements
will lead to products for non-minimal length coset representatives or
non-reduced expressions.

\begin{corollary}\label{c:bruhat.covers.B}
If $v <w$ in Bruhat order on $\minreps$ and $l(w)=l(v)+1$, then one of
the following must hold:
\begin{enumerate}
\item The partition $\pi(v)$ is obtained from $\pi(w)$ by removing a
single outer corner square so that the remaining shape is in
$\partitions$.

\item The partition $\pi^{-1}(v)$ is obtained from $\pi^{-1}(w)$ by
removing a single inside square $s$ in row $i$ of the Ferrers diagram
so that the generators corresponding to the squares to the right of
square $s$ all commute up to join onto the end of an existing row
above row $i$, and furthermore, the resulting product of segments can
be put into canonical form using the commutation relations in
Lemma~\ref{l:commutation.B} and Lemma~\ref{l:C.commutations.B} in such
a way that they correspond with a partition in $\partitions$
with one fewer square.
\end{enumerate}
\end{corollary}

\begin{lemma}\label{l:youngs.lattice.B}
If $\lambda \in \partitions(B_{n})$ and $\lambda_{1}<n$, then the
interval in Bruhat order between $\id$ and $\pi ^{-1} (\lambda)$ is
isomorphic as posets to the interval from the empty partition to
$\lambda $ in Young's lattice on strict partitions.
\end{lemma}

\begin{proof}
Since $\lambda \in \partitions(B_{n})$ and $\lambda_{1}<n$, $\lambda$
must be a strict partition contained in the staircase
$(n-1,\dotsc,2,1)$.  Note that all strict partitions whose parts are
of length less than $n$ occur in $\partitions(B_{n})$ and by
Corollary~\ref{c:bruhat.shape.covers} all the covering relations from
Young's lattice restricted to this set correspond with covering
relations in Bruhat order on $\minreps$.  Therefore, $w=\pi ^{-1}(\lambda
) \geq v= \pi ^{-1}(\mu )$ for all $\mu \in \partitions(B_{n})$ and
$\mu \subset \lambda$.  

In order to prove the converse, it suffices to show that deleting any
generator from $r(w)=\cdots
\segment{}{0}(\lambda_{3})\segment{}{1}(\lambda_{2})\segment{}{0}(\lambda_{1})$
which doesn't correspond with a $\partitions$-removable corner leaves a reduced
expression for an element of $\widetilde{W}$ that is not a minimal
length coset representative in $\minreps$ of length $l(w)$.  Then by
induction on $l(w)$ the result follows.

Identify $r(w)$ with the filling of $\lambda$ as described in
Example~\ref{ex:labeling.B}. Say a generator in column $i$ and row $j$
is deleted from $\lambda$ using matrix notation.

If $(i,j)$ corresponds with a corner of the partition, then either
that corner is $\partitions$-removable or $\lambda_{j} = \lambda_{j+1}+1$.  If
$\lambda_{j} = \lambda_{j+1} +1$ then deleting $s_{i}$ leaves two
equal parts of the remaining partition.  By
Lemma~\ref{l:C.commutations.B}, the product of two equal length
segments can be rewritten as the product of two segments times $s_{1}$
or $s_{0}$ on the right.  By Lemma~\ref{l:commutation.B}, the $s_{0}$
or $s_{1}$ commutes up in $\lambda$ increasing its index as it passes
each row, but since $\lambda$ is a strict partition with each
$\lambda_{k}<n$, it can never glue onto the end of an existing segment
above.  Therefore, removing a corner that is not $\partitions$-removable leaves
an expression for a non-minimal length element in its coset.

Similarly, if $(i,j)$ corresponds with an internal square of $\lambda
\subset (n-1,\dotsc,2,1)$, then we claim $r(w)$ with the generator in
$(i,j)$ removed has a right factor which is not an element in
$\minreps$.  To observe the right factor, note that after deleting
$(i,j)$, the adjacent $s_{i+1}$ in row $j$ commutes up to the next row
in $\lambda$.  Applying Lemma~\ref{l:commutation.B} and the fact that
$i+1\leq \lambda_{j}<\lambda_{j-1}\dotsb <\lambda_{1}<n$, one has
\[
s_{i+1} \segment{}{0}(\lambda_{j-1})\ldots
\segment{}{0}(\lambda_{3})\segment{}{1}(\lambda_{2})\segment{}{0}(\lambda_{1})
=  \segment{}{0}(\lambda_{j-1})\ldots
\segment{}{0}(\lambda_{3})\segment{}{1}(\lambda_{2})\segment{}{0}(\lambda_{1}) s_{i+j}
\]
with $1<i+j\leq n$.  Thus, $r(w)$ with the generator in $(i,j)$
removed does not correspond with a covering relation in Bruhat order.
\end{proof}

\begin{lemma}\label{l:chains.B}
If $\lambda \in \partitions(B_{n})$ is a single row ($\lambda = (j)$),
then the interval in Bruhat order between $\id$ and $\pi ^{-1}
(\lambda)$ is isomorphic as posets to the interval from the empty
partition to $\lambda $ in Young's lattice.  In particular, this
interval is a totally ordered chain.
\end{lemma}

\begin{proof}
If $\lambda$ is a single row of length $j$, then $r(\pi ^{-1}(\lambda
))= \segment{}{0}(j)$. Observe directly from the relations on the
generators that deleting any but the leftmost generator from
$\segment{}{0}(j)$ either leaves a non-reduced expression or a
non-minimal length coset representative by definition.
\end{proof}

Consider the Bruhat order on $\minreps$ expressed in terms of
partitions in $\partitions(B_{n})$ for $n\geq 3$.  By
Lemma~\ref{l:youngs.lattice.B}, the smallest 4 ranks of the Bruhat
Hasse diagram are determined by Young's lattice on strict partitions:

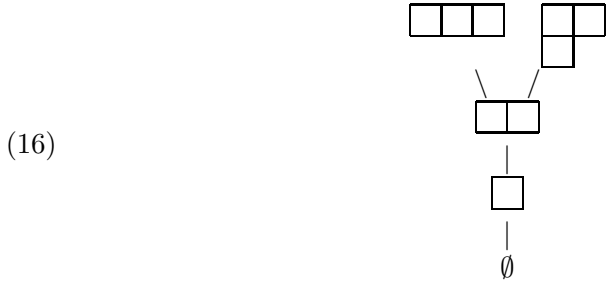
\begin{equation}\label{e:hasse.diag.B}
\begin{array}{c}
\tableau{{} & {} &{}}   \hspace{.2in} \tableau{{} & {}\\ {} & } \\ 
\backslash  \hspace{.2in} /
\\
\tableau{{} & {}}\\
| \\
\tableau{{} }\\
| \\
\emptyset
\end{array}
\end{equation}

\begin{theorem}\label{t:rationally.smooth.B} Let $\widetilde{W}$ be
the affine Weyl group of type $B_{n}$ and let $w \in \minreps$.  For
$n\geq 4$, $w$ is palindromic if and only if $\pi(w)$ is a one row
shape $(j)$ with $0\leq j<2n$ or a staircase shape $(k, k-1,
k-2,\dots, 1)$ for $1<k< n$.  For $n=3$, $w$ is palindromic if and
only if $\pi(w)$ is a one row shape or the staircase shape $(2,1)$ or
the square $(3,3,3)$.
\end{theorem}

\begin{rem}
The number of palindromic elements in $\minreps$ for type $B_{n}$ is
therefore 8 for $n=3$ and $3n-2$ for $n\geq 4$.  The one row shapes
with $n\leq j\leq 2n-2$ and the $(3,3,3)$ shape in $n=3$ correspond
with rationally smooth Schubert varieties which are not actually
smooth.  All others are closed parabolic orbits, hence they are smooth
by Theorem~\ref{t:billey.mitchell}. 
\end{rem}

\begin{proof}
By definition, $w$ is palindromic if and only if its Poincar\'e
polynomial $P_{w}(t)=\sum_{v \leq w} t^{l(v)}$ is palindromic.  By
Lemma~\ref{l:chains.B}, we know that $P_{w}(t) = 1 +t +\dotsb +t^{j}$
if $\pi (w)$ is a one row shape, which is clearly palindromic.  By
Lemma~\ref{l:youngs.lattice.B}, if $\pi (w)$ is a staircase shape $(k,
k-1, k-2,\dots, 1)$ for $1<k< n$ and $n\geq 3$, then
\[
P_{w}(t) = \prod_{j=1}^{k} (1+t^{j})
\]
which is easily seen to be palindromic by induction.  In the case
$n=3$, one can verify that for $w=s_{3}s_{2}s_{0}\cdot
s_{3}s_{2}s_{1}\cdot s_{3}s_{2}s_{0}$ corresponding with
$\pi(w)=(3,3,3)$, the Poincar\'e polynomial for $X_{w}$ is
\begin{align*}
P_{w}(t)= &
t^9 + t^8 + t^7 + 2t^6 + 2t^5 + 2t^4 + 2t^3 + t^2 + t^1 + 1\\
= & (1+t)(1+t^{2})(1+t^{3}+t^{6}) + t^{4} + t^{5}
\end{align*}
which is palindromic.\footnote{It is interesting to observe that the
interval below $w=\pi ^{-1}{(3,3,3)}$ contains all the elements
corresponding with the strict partitions contained in $(3,3,3)$ plus
the two additional elements $\pi ^{-1} (5)$ and $\pi ^{-1} (4)$ in
Bruhat order which contribute symmetrically to $P_{w}(t)$.}

Conversely, we will show that if $\pi (w)$ is not a one row shape or a
staircase shape  or   $(3,3,3)$  in  the case  $n=3$,  then  for  some
$0<j<l(w)$   one    pair  of  the    coefficients   in   
\[
P_{w}(t)   =
\sum_{i=0}^{l(w)} a_{i}t^{i}
\]
has the property that $a_{j} \neq a_{l(w)-j}$.  We can calculate all
$a_{j}$ for $0<j\leq n$ by computing the number of strict partitions
that fit inside the Ferrers diagram for $\pi (w)$ since there are no
repeated parts in any partition in $\partitions(B_{n})$ of size $n$ or
less and $\partitions(B_{n})$ contains all such strict partitions.
Therefore, in order to show $P_{w}(t)$ is not palindromic, it suffices
to show there exists some $ 1<j \leq n$ such that $a_{l(w)-j}$ is
bigger than the number of strict partitions below $\pi (w)$ of size
$j$.

Assume $\lambda = \pi(w)$ has at least two rows.  From
\eqref{e:hasse.diag.B}, we know the branching number $b_{W}=2$ so
$a_{1}=a_{2}=1$.  If $\lambda$ has two or more $\partitions$-removable outside
corners, then $a_{l(w)-1} \geq 2$ so $w$ is not palindromic.
Similarly, if $\lambda$ has one $\partitions$-removable corner and removing it
leaves a shape with two $\partitions$-removable corners, then $a_{l(w)-2} \geq 1$ so
$w$ is not palindromic.  One of these two cases occurs for all
partitions $\lambda \in \partitions(B_{n})$ with at least two rows and
$\lambda_{1}>n$ or if there exists an $i$ such that $\lambda_{i}
-\lambda_{i+1}\geq 2$.  Therefore, it remains to show that if $\lambda
= (n,n,\ldots, n)=(n^{k})$ or $\lambda = (n^{k},n-1,\ldots,1)$ for
some $k>1$ then $P_{w}(t)$ is not palindromic excluding the case
$\lambda =(3,3,3)$ and $n=3$.

If $\lambda = (n^{k})$ then the interval below $\lambda$ in Young's
lattice restricted to $\partitions(B_{n})$ is self dual, and hence
rank symmetric, since we can pair any $\gamma \subset \lambda $ with
its compliment inside the rectangle of size $n \times k$.  Therefore,
for each $k>1$, in order to show $w$ is not palindromic, we simply
need to find a single partition $\mu \in \partitions(B_{n})$ such that
$\pi ^{-1}(\mu) < w$ in Bruhat order, $\mu \not \subset \lambda$ in
Young's lattice, and $\ell(\mu ) \geq l(w)-n$.

Assume first that $1<k<n$. Let $i=n-k$.  Then,
$w>\pi^{-1}(n^{k-1},n-1)$ and deleting the generator in column $i$ on
row $k$ in the reduced expression corresponding with $(n^{k-1},n-1)$
gives an expression where the leftmost $s_{i+1}$ commutes right past
$k-2$ segments of length $n$, increasing its index as it commutes past
each segment, to become an $s_{i+1+k-2}=s_{n-1}$ which glues onto the
rightmost $\segment{}{0}(n)$.  Similarly, the leftmost $s_{i+2}$
commutes right pas $k-3$ segments and glues onto the rightmost
$\segment{}{1}(n)$ to become $\segment{}{1}(n+1)$, etc for the left
most $s_{i+3},\ldots s_{n-1}$.  The remaining expression corresponds
with the partition $\mu =((n+1)^{k-1},n-k-1) \not \subset (n^{k})$,
but $v=\pi ^{-1}(\mu)<w$ in Bruhat order and $\ell(v)=l(w)-2$ so $w$
is not palindromic.  For example, if $n=6, k=4$, then if $w=\pi
^{-1}(6,6,6,6)$, we have $\mu = (7,7,7,1)$ and $v=\pi ^{-1}(\mu)$ is
obtained from $\pi^{-1}(6,6,6,5)$ by deleting one more generator and
applying commutation relations:
\[\lambda =
\tableau{
0 & 2 & 3 & 4 & 5 & 6\\
1 & 2 & 3 & 4 & 5 & 6\\
0 & 2 & 3 & 4 & 5 & 6\\
1 & 2 & 3 & 4 & 5 & 6 \\
}
\text{\hspace{.1in} $>$ \hspace{.1in} }
\tableau{
0 & 2 & 3 & 4 & 5 & 6\\
1 & 2 & 3 & 4 & 5 & 6\\
0 & 2 & 3 & 4 & 5 & 6\\
1 &  & 3 &  4&  5 & \\
} = 
\tableau{
0 & 2 & 3 & 4 & 5 & 6 & 5\\
1 & 2 & 3 & 4 & 5 & 6 & 5\\
0 & 2 & 3 & 4 & 5 & 6 & 5\\
1 &  &  &  &  & \\
} = \mu 
\]


Assume next that $\lambda =(n^{n})$, we claim that $\mu =
((n+2)^{n-2})$ is a partition in $\partitions(B_{n})$ which is not
contained in $\lambda$ such that $\pi ^{-1}(\mu ) <w$ in Bruhat order
and $l(\pi ^{-1}(\mu))=l(w)-4$.  For $n\geq 4$, this implies $a_{4}<
a_{l(w)-4}$, so $w$ is not palindromic.  To prove the claim,
observe that by Lemma~\ref{l:youngs.lattice.B}, $w>\pi
^{-1}(n^{n-2},n-1,n-2)$.  Let $v'$ be the element obtained by deleting
the generator corresponding to row $n-1$ and column 1 in $\pi
^{-1}(n^{n-2},n-1,n-2)$.  Assume $n$ is even for the sake of notation
(the case $n$ odd is the same except for a 0,1 switch).  Then, by
definition
\[
v' = s_{n-2}\cdots s_{2}s_{1} \cdot s_{n-1} \dotsb s_{3}s_{2} \cdot
\pi ^{-1}(n^{n-2}).  
\]
Note that by Lemma~\ref{l:commutation.B}, $s_{2} \cdot \pi
^{-1}(n^{n-2}) = \pi ^{-1}(n+1, n^{n-3})$ since the index on the
generator $s_{2}$ increases by one each time it commutes with a
$\segment{}{0}(n)$ or $\segment{}{1}(n)$, after commuting with $n-3$
segments $s_{2}$ becomes $s_{n-1}$ and $s_{n-1} \cdot \segment{}{0}(n)
= \segment{}{0}(n+1)$.  After moving the $s_{2}$, the last $s_{1}$
commutes right as far as possible and becomes an $s_{n-2}$ which glues
onto the rightmost segment.  Continuing in the same way we can commute
all the generators remaining in the expression for $v'$ above in the
last two rows of $(n^{n})$ up to glue onto rows $1,\ldots ,n-2$ adding
two additional generators to each row, so $\pi ^{-1}(\mu)=v <w$.

Similarly, if $\pi (w)$ is any partition in $\partitions(B_{n})$ with
$n\geq 4$ of the form $=(n^{k},\lambda_{k+1},\dotsc ,\lambda_{m})$
with $k>n$, then $r(w)$ has a right factor of $\pi ^{-1} (n^{n})$ so
as above $w$ is bigger than $v=\pi ^{-1}((n+2)^{n-2},
n^{k-n},\lambda_{k+1},\dotsc ,\lambda_{m})$ and $l(v)=l(w)-4$, so
$w$ is not palindromic.

If $n=3$, we claim any affine partition $(3^{k},\lambda_{k+1},\dotsc
,\lambda_{m})$ is not palindromic for $k> 3$.  To see this note that
$(3^{5})$ covers $(4,4,3,3)$, so $(3^{k},\lambda_{k+1},\dotsc
,\lambda_{m})$ covers at least 2 elements in Bruhat order, hence not
non-palindromic, for $k\geq 5$ by an argument similar to the proof of
Lemma~\ref{l:exceptional.finite}.  For $k=4$, it is easy to verify the
claim on the 4 possible affine partitions.  For example, $(3^{4})$
covers $(3,3,3,2)$ and $(4,3,3,1)$ chopping out the 3 on the 3rd
row. Similarly, $(3^{4},2)$ covers $(4,3,3,3)$ and $(4,4,3,2)$.

Finally, we need to address the case $\lambda
=(n^{k},n-1,n-2,\ldots,2,1)$ for $1\leq k<n$ and $n\geq 4$.  We claim
$w$ covers two elements of length $l(w)-1$ in $\minreps$.  The first
element is obtained by deleting the last row of $\lambda$.  The second
element, indexed by the partition $((n+1)^{k-1},
n,n-1,\dotsc,\widehat{k},\dotsc,2,1)$, is obtained from $\lambda$ by
deleting the corner square in row $n$ column $k$ to obtain the
expression (suppressing the 0,1 subscripts so we don't need to
consider the parity of $n$):
\[
\segment{}{}(1) \segment{}{}(2)\dotsb \segment{}{}(k-1)\segment{}{}(k-1)
\cdot \pi ^{-1}(n^{k},n-1,n-2,\ldots,k+1)
\]
Now, by Lemma~\ref{l:C.commutations.B},
$\segment{}{1}(k-1)\segment{}{0}(k-1)=
\segment{}{1}(k-2)\segment{}{0}(k-1)\cdot s_{1}$ and moving $s_{1}$ up
past $n-2$ rows becomes $s_{n-1}$ by Lemma~\ref{l:commutation.B} which
glues onto row 1 to become a $\segment{}{0}(n+1)$.  Swapping $s_{1}$
with $s_{0}$ results in the same partition shapes.  After rectifying
the two segments of length $k-1$ there are two segments of length
$k-2$ which rectifies by adding an additional generator to row 2,
etc.  Applying the same technique repeatedly we can rectify the whole
expression to obtain an equivalent reduced expression corresponding to
the shape obtained from $\lambda$ by removing the corner square from
each row from $n$ to row $n+k-1$ and adding one square to rows 1
through $k-1$.
\end{proof}

Recall the definition of $m_{W}$ from Section~\ref{s:palindromics}.

\begin{corollary}\label{c:mw.B}
For $W$ of type $B_{n}$, we have $m_{W}=4$ for $n\geq 4$ and  $m_{W}=2$ for $n=3$.  
\end{corollary}

\section{Type $C$}\label{s:typeC}

Let $\partitions(C_{n})$ be the set of partitions whose parts are
bounded by $2n$ and the parts of length less than or equal to $n$ are
strictly decreasing.  The generating function for such partitions is
\[
G_{\Cn}(x)=\frac{(1+x)(1+x^{2})\cdots (1+x^{n})}{(1-x^{n+1}) (1-x^{n+2}) \cdots (1-x^{2n})}
\]

\begin{lemma}\label{l:Cpartitions}
We have the following generating function identity with Bott's formula
from \eqref{e:bott.formula}:
\[
G_{\Cn}(x) = \frac{1}{(1-x)(1-x^{3})(1-x^{5})\cdots (1-x^{2n-1})}.
\]
\end{lemma}

\begin{proof}
Apply induction on $n$. 
\end{proof}

The Dynkin diagram of $\widetilde{C}_{n}$ with the extra generator labeled
$s_{0}$ is adjacent to $s_{1}$ so $\minreps$ will have Type I segments
and fragments as described in Section~\ref{s:canonical}.  The
parabolic subgroup generated by $S=\{s_{1},s_{2}, s_{3}\dotsc , s_{n}
\}$ corresponds with the finite Weyl group of type $C_n$.  Let
$J=\{s_{2},\dotsc , s_{n} \}$.  Then the minimal length coset
representatives for $W/W_{J}$ are the fragments
\begin{equation}\label{e:minrepsC}
\begin{array}{cccc}
\id & \\
 s_{i} \dotsb s_{3}s_{2}s_{1}  &	\textrm{for} \ 1\leq i \leq n\\
s_{i} \dotsb s_{n-1}s_{n} s_{n-1}\dotsb s_{3}s_{2} s_{1} &\textrm{for} \ 1\leq i < n
\end{array}
\end{equation}
The fact that these are all the minimum length coset representatives
for $W_{I}/W_{J}$ is easily checked using the signed permutation
description of the type $C_n$ Weyl group elements \cite{b-b}.   
We will use the reduced expressions in \eqref{e:minrepsC} to define
the \textit{segments} $\segment{i}{}(j)$.  Note, in this case, each
word in \eqref{e:minrepsC} has unique length, so we will denote
$\segment{1}{}(j)$ by simply $\segment{}{}(j)$.  Set
$\segment{}{}(0)=\id$, and for $1\leq j \leq 2n$ set
\begin{equation}\label{e:Csegments.1}
\segment{}{}(j) = \begin{cases}
s_{j-1} \dots  s_{2} s_{1} s_{0} &	 1\leq j\leq n+1\\
s_{2n-j+1} \dotsb s_{n-1}s_{n} s_{n-1}\dotsb s_{2}  s_{1}  s_{0}&
n+1<j\leq 2n.
\end{cases}
\end{equation}

The proofs of the next two lemmas are analogous to
Lemmas~\ref{l:commutation.B} and~\ref{l:C.commutations.B} in type $B$.

\begin{lemma}\label{l:commutation.C}
We have the following commutation rules for $s_{i} \cdot \segment{}{}(j)$ for all
$1 \leq i \leq n $ and $1\leq j\leq 2n$:
\begin{equation}\label{e:Csegments}
s_{i} \cdot \segment{}{}(j) = \begin{cases}
\segment{}{}(j) \cdot s_{i} &	 1\leq j < i  \text{ or  }   2n-i+1<j \leq 2n\\
\segment{}{}(j+1) &	 j = i \text{ or }j = 2n -i\\
\segment{}{}(j-1) &   j = i+1  \text{ or }j = 2n -i +1\\
\segment{}{}(j) \cdot s_{i+1} &  i +1 < j < 2n-i.
\end{cases}
\end{equation}
\end{lemma}

\begin{lemma}\label{l:C.commutations.C}
We have the following product rules for segments in type $C_{n}$:

\begin{enumerate}
\item For $1\leq j\leq n$,  
\[
\segment{}{}(j) \cdot \segment{}{}(j) =
\segment{}{}(j-1)\cdot  \segment{}{}(j)\cdot s_{1}.
\]

\item For $n+1\leq j \leq 2n$, 
\[
\segment{}{}(j+1) \cdot  \segment{}{}(j) =
\segment{}{}(j)\cdot \segment{}{}(j)\cdot s_{1}.
\]
\end{enumerate}
\end{lemma}





\begin{theorem}\label{t:surjective.Cn}
In type $C_{n}$, $\pi : \minreps \longrightarrow \partitions =
\partitions(\Cn)$ is a length preserving bijection.
\end{theorem}

\begin{proof} 
As in the proof of Theorem~\ref{t:surjective.Bn}, we only need to show
that the allowed pairs in type $\Cn$ have parts bounded by $2n$ and
are strictly increasing if the larger part has length less than or
equal to $n$.  This follows directly from
Lemma~\ref{l:C.commutations.C}.
\end{proof}

\begin{corollary}\label{c:allowed.pairs.C}
The product $\segment{}{}(j) \cdot \segment{}{}(k) \in \minreps$ and
$\ell(\segment{}{1}(j) \cdot \segment{}{0}(k))=j+k$ if and only if
$1\leq j<k\leq 2n $ or $n< j=k \leq 2n$.  
\end{corollary}

\begin{ex}\label{ex:labeling.C}
If $n=4$, the partition $(7,5,5,3,1) \in \partitions(C_{n})$
corresponds with $s_0\cdot s_{2}s_{1}s_{0} \cdot
s_{4}s_{3}s_{2}s_{1}s_{0} \cdot s_{4}s_{3}s_{2}s_{1}s_{0} \cdot
s_{2}s_{3}s_{4}s_{3}s_{2}s_{1}s_{0}$.  Pictorially, we have 
\[
\tableau{
0 & 1 & 2 & 3 & 4 & 3 & 2 & \\
0 & 1 & 2 & 3 & 4 &  &  & \\
0 & 1 & 2 & 3 & 4 &  &  & \\
0 &  1 &2 &  &  &  &  & \\
0 &  &  &  &  &  & }
\]
\end{ex}
\noindent where again the corresponding reduced expression is read
along the rows from right to left, bottom to top.  The columns contain
their column number minus one up to column $n+1$ and for $n+1< j \leq
2n-1$ column $j$ contains $2n-j+1$.

Bruhat order on $\minreps$ in type $C$ is very similar to type $B$.  The
proofs for the next two lemmas are very similar to the proofs of
Lemmas~\ref{l:youngs.lattice.B} and~\ref{l:chains.B}.  

\begin{lemma}\label{l:youngs.lattice.C}
If $\lambda \in \partitions(C_{n})$ and every $\lambda_{k} \leq n$,
then the interval in Bruhat order between $\id$ and $\pi ^{-1}
(\lambda)$ is isomorphic as posets to the interval from the empty
partition to $\lambda $ in Young's lattice on strict partitions.
\end{lemma}

\begin{lemma}\label{l:chains.C}
If $\lambda \in \partitions(C_{n})$ is a single row, then the interval
in Bruhat order between $\id$ and $\pi ^{-1} (\lambda)$ is isomorphic
as posets to the interval from the empty partition to $\lambda $ in
Young's lattice.  In particular, this interval is
a totally ordered chain.
\end{lemma}

Consider the elements in $\minreps$ determined by partitions in
$\partitions(C_{n})$ for $n\geq 2$.  By
Lemma~\ref{l:youngs.lattice.C}, the smallest 5 ranks of Bruhat order
on $\minreps$ is the same as Young's lattice on strict partitions of size
at most 4

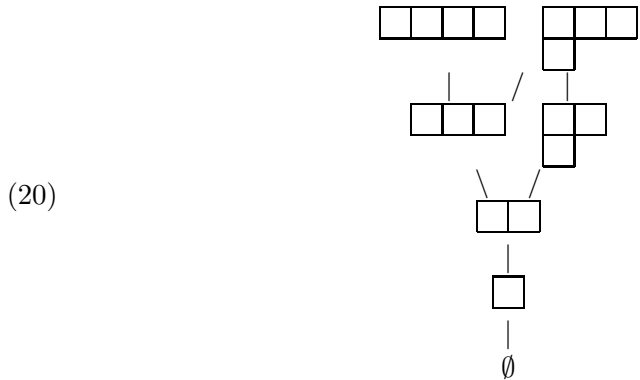
\begin{equation}\label{e:hasse.diag.C}
\begin{array}{c}
\tableau{{} & {} &{} & {}}   \hspace{.2in} \tableau{{} & {} & {}\\ {} & } 
\\ 
|  \hspace{.3in} /  \hspace{.2in} |
\\
\tableau{{} & {} &{}}   \hspace{.2in} \tableau{{} & {}\\ {} & } \\ 
\backslash  \hspace{.2in} /
\\
\tableau{{} & {}}\\
| \\
\tableau{{} }\\
| \\
\emptyset
\end{array}
\end{equation}

\begin{theorem}\label{t:rationally.smooth.C} Let $W$ be
the Weyl group of type $C_{n}$ for $n\geq 2$ and let $w \in
\minreps$.  Then $w$ is palindromic if and only if $\pi(w)$ is a one
row shape $(j)$ for $0\leq j\leq 2n$ or a staircase shape $(k, k-1,
k-2,\dots, 1)$ for $1<k\leq n$.
\end{theorem}

\begin{rem}
The number of palindromic elements in $\minreps$ for type $C_{n}$ is
therefore $3n$ for $n\geq 2$.  The Schubert varieties indexed by one
row shapes with $2\leq j\leq 2n$ are rationally smooth but not
actually smooth.  All of the staircase shapes are closed parabolic
orbits, hence they are smooth by Theorem~\ref{t:billey.mitchell}. 
\end{rem}

\begin{proof}
The proof of Theorem~\ref{t:rationally.smooth.C} is similar to the
proof of Theorem~\ref{t:rationally.smooth.B}. By
Lemma~\ref{l:chains.C}, we know that $P_{w}(t) = 1 +t +\dotsb +t^{j}$
if $\pi (w)$ is a one row shape, which is clearly palindromic.  By
Lemma~\ref{l:youngs.lattice.C}, if $\pi (w)$ is a staircase shape $(k,
k-1, k-2,\dots, 1)$ for $1<k\leq n$ then
\[
P_{w}(t) = \prod_{j=1}^{k} (1+t^{j})
\]
which is palindromic.

Conversely, assume $\lambda = \pi(w)$ has at least two rows and is not
a staircase shape for $k\leq n$.  From \eqref{e:hasse.diag.C}, we know
$a_{1}=a_{2}=1$.  If $\lambda$ has two or more $\partitions$-removable
outside corners, then $a_{l(w)-1} \geq 2$ so $w$ is not palindromic.
Similarly, if $\lambda$ has one $\partitions$-removable corner and
removing it leaves a shape with two $\partitions$-removable corners,
then $a_{l(w)-2} \geq 1$ so $w$ is not palindromic.  One of these two
cases occurs for all partitions $\lambda \in \partitions(C_{n})$ with
at least two rows and $\lambda_{1}>n+1$ or if there exists an $i$ such
that $\lambda_{i} -\lambda_{i+1}\geq 2$.  Therefore, it remains to
show that if $\lambda = (n+1)^{k}$ for $k\geq 2$ or $\lambda =
((n+1)^{k},n,n-1,\ldots,1)$ for some $k\geq 1$ then $P_{w}(t)$ is not
palindromic.  Note, these are the only two families of thin elements
which are not palindromic in type $C$.

If $\lambda = ((n+1)^{k})$ for $k\leq n$ then $\lambda$ covers only
one element $\mu = ((n+1)^{k-1}, n)$ in Bruhat order.  For $k<n$,
$\mu$ covers both $((n+1)^{k}, n-1)$ and $((n+2)^{k-1}, n-k)$.  The
later covering relations is seen by deleting the $s_{n-k}$ from the
$\segment{}{}(n)$ factor which allows the $s_{n-k+1}, \ldots, s_{n-1}$
in the bottom row to commute up adding one generator to each other
factor.  For $k=n$, $\mu$ covers $((n+2)^{n-1})$ by deleting the
$s_{0}$ in the bottom row.  Therefore, for all $k\geq 2$ there are at
least two elements of rank $|\lambda |-2$ below $\lambda =
((n+1)^{k})$ in Bruhat order.

Assume now $\lambda = ((n+1)^{k},n,n-1,\ldots,1)$.  If $k\geq n$ then
it follows from the previous case that $\lambda$ is not
palindromic.  If $0<k<n$, then $\lambda$ covers two elements
$((n+1)^{k},n,n-1,\ldots,2)$ and
$((n+2)^{k-1},n,n-1,\ldots,k+2,k,k-1,\ldots 1)$.  The later element is
obtained from $\lambda$ by knocking out the $s_{0}$ in the
$\segment{}{}(k+1)$ factor and allowing the $s_{1},\dots, s_{k}$
generators to commute up adding one generator to each of the
$\segment{}{}(n+1)$ factors.  Hence, for all $k\geq 1$, $\lambda $ is
not palindromic.



\end{proof}

\begin{corollary}\label{c:mw.C}
For $W$ of type $C_{n}$, we have $m_{W}=2$.  
\end{corollary}

\section{Type $D$}\label{s:typeD}

Let $\partitions(D_{n})$ be the set of partitions whose parts are
bounded by $2n-2$, the parts of length at most $n-2$ are
strictly decreasing, and the parts of length $n-1$ are two colored,
but at most one color can appear in any partition in
$\partitions(D_{n})$.  The generating function for such partitions is
\[
G_{\Dn}(x)=\frac{(1+x)(1+x^{2})\cdots (1+x^{n-2})(1+x^{n-1})}
{(1-x^{n-1})(1-x^{n}) (1-x^{n+1})\cdots (1-x^{2n-2})} 
\]

\begin{lemma}\label{l:Dpartitions}
We have the following generating function identity with Bott's formula
from \eqref{e:bott.formula}:
\[
G_{\Dn}(x) = \frac{1}{(1-x^{n-1})(1-x)(1-x^{3})(1-x^{5})\cdots (1-x^{2n-3})}.
\]
\end{lemma}

\begin{proof}
The formula holds for $n=4$ and the general case follows by induction.
\end{proof}

The Coxeter graph for $\widetilde{D}_{n}$ has an involution
interchanging $s_{0}$ and $s_{1}$ so this is a Type II Coxeter graph.
The segments described in Section~\ref{s:canonical} for type $D$ are
given as follows: for $1\leq j \leq 2n-2$, $z\in \{b,c\}=\{\mathrm{blue,
crimson}\}$ set
\begin{equation}\label{e:Dsegments.1}
\segmentD{z}{1}(j) = \begin{cases}
s_{j} \dots  s_{3}s_{2} s_{1} &	 1\leq j\leq n-2\\
s_{n-1} s_{n-2} \dots  s_{3}s_{2} s_{1} &  j=n-1 \text{ and  } z=b\\
s_{n} s_{n-2} \dots  s_{3}s_{2} s_{1} &  j=n-1 \text{ and  } z=c\\
s_{2n-j-1} \dotsb s_{n-2}s_{n-1}s_{n} s_{n-2}\dotsb s_{4}s_{3}s_{2}  s_{1} &
n\leq j \leq  2n-2.
\end{cases}
\end{equation}
Note, there are two distinct 1-segments of length $n-1$, namely
$\segmentD{b}{1}(n-1) $ and $\segmentD{c}{1}(n-1)$.  The superscript
$z=b,c$ distinguishes these two cases.  Otherwise, for $j\neq n-1$ the
$z$ is irrelevant and can be omitted if it is convenient.

Similarly, define $\segmentD{z}{0}(j)$ for all $1\leq j\leq 2n-2$ and
$z\in \{b,c\}$ to be the reduced expression obtained from
$\segmentD{z}{1}(j)$ by interchanging $s_{0}\leftrightarrow s_{1}$ and
$s_{n-1}\leftrightarrow s_{n}$.  For example,
$\segmentD{b}{1}(n-1)=s_{n-1} s_{n-2} \dots s_{3}s_{2} s_{1}$ so
$\segmentD{b}{0}(n-1)=s_{n} s_{n-2} \dots s_{3}s_{2} s_{0}$.  Note the
longest 0-segment $\segmentD{z}{0}(2n-2)=s_{0}s_{2} \dotsb
s_{n-2}s_{n-1}s_{n} s_{n-2}\dotsb s_{4}s_{3}s_{2} s_{0}$ both starts
and ends in $s_{0}$.

\begin{lemma}\label{l:commutation.D}
In type $\widetilde{D}_{n}$, we have the following commutation rules
for $s_{i} \cdot \segmentD{z}{1}(j)$ for all $1\leq j\leq 2n-2$ and $z
\in \{b,c\}$.  For $1 \leq i \leq n-2$,
\begin{align}\label{e:Dsegments}
&s_{i} \cdot \segmentD{z}{1}(j) = \begin{cases}
\segmentD{z}{1}(j) \cdot s_{i} &	 1\leq j < i-1   \text{ or  }   2n-i-1<j \leq 2n-2\\
\segmentD{z}{1}(j+1) &	 j = i-1 \text{ or  } j = 2n -i-2\\
\segmentD{z}{1}(j-1) &   j = i \text{ or  } j=2n-i-1\\
\segmentD{z}{1}(j)\cdot s_{i+1} &  i <j  <  2n-i-2
\end{cases}
\end{align}
For $i=n-1$, we have
\begin{align}\label{e:Dsegments.n-1}
&s_{n-1} \cdot \segmentD{z}{1}(j) = \begin{cases}
\segmentD{}{1}(j) \cdot s_{n-1} &	 1\leq j < n-2  \\
\segmentD{b}{1}(j+1) &	 j = n-2 \\
\segmentD{}{1}(j+1) &	 j = n-1 \text{ and  } z=c\\
\segmentD{}{1}(j-1) &   j = n-1   \text{ and  } z=b\\
\segmentD{c}{1}(j-1) &   j = n  \\
\segmentD{}{1}(j) \cdot s_{n} &  n < j \leq 2n-2.
\end{cases}
\end{align}
For $i=n$, the commutation relations for $s_{n} \cdot
\segmentD{z}{1}(j)$ are similar to those in \eqref{e:Dsegments.n-1}
interchanging $b \leftrightarrow c$ and $s_{n} \leftrightarrow s_{n-1}$.  In
particular, for $n<j\leq 2n-2$,
\[
s_{n} \cdot \segmentD{z}{1}(j) = \segmentD{z}{1}(j) \cdot s_{n-1}.
\]
 For $i=0$
\[
\segmentD{z}{0}(1) \cdot \segmentD{z}{1}(1) =
s_{0}s_{1} = s_{1}s_{0} = \segmentD{z}{1}(1) \cdot \segmentD{z}{0}(1).
\]
Similarly, commutation rules for $s_{i} \cdot \segmentD{z}{0}(j)$ for
all $0 \leq i \leq n $, $z \in \{b,c\}$ and $1\leq j\leq 2n-2$ are
obtained from \eqref{e:Dsegments} and \eqref{e:Dsegments.n-1} by
interchanging $s_{0}\leftrightarrow s_{1}$ and $s_{n-1}\leftrightarrow
s_{n}$. 
\end{lemma}

\begin{proof}
These follow directly from the commutation relations among the
generators determined by the Dynkin diagram.
\end{proof}


\begin{lemma}\label{l:C.commutations.D}
We have the following product rules for segments in type
$\widetilde{D}_{n}$:

\begin{enumerate}
\item For $1\leq j <n-1$, we have
\[
\segmentD{}{1}(j) \cdot \segmentD{}{0}(j) =
\segmentD{}{1}(j-1)\segmentD{}{0}(j)\cdot s_{1}.
\]

\item For $j=n-1$, we have 
\[
\segmentD{b}{1}(n-1) \cdot \segmentD{c}{0}(n-1) =
\segmentD{}{1}(n-2)\segmentD{c}{0}(n-1)\cdot s_{1}
\]
and 
\[
\segmentD{c}{1}(n-1) \cdot \segmentD{b}{0}(n-1) =
\segmentD{}{1}(n-2)\segmentD{b}{0}(n-1)\cdot s_{1}.
\]
\item For $n\leq j <2n-1$, we have 
\[
\segmentD{}{1}(j+1) \cdot \segmentD{}{0}(j) =
\segmentD{}{1}(j)\segmentD{}{0}(j)\cdot s_{1}.
\]
\end{enumerate}
Similar rules hold interchanging $0 \leftrightarrow 1$ above.
\end{lemma}

\begin{proof}
If $1\leq j < n-1$, then the segments $\segmentD{}{0}(j)$ and
$\segmentD{}{1}(j)$ involve only a proper subset of the generators of
$\widetilde{S}$.  This proper subset generates a parabolic subgroup
isomorphic to a parabolic subgroup of the affine Weyl group of type
$\widetilde{B}_{n}$.  Therefore, the stated relation holds by
Lemma~\ref{l:C.commutations.B}.

Assume $n\leq j<2n-1$ and let $t=2n-j-1$.  Then, 
\begin{align*}
\segmentD{}{0}(j)&=s_{t}s_{t+1}\dotsb s_{n-2}s_{n-1}s_{n}s_{n-2}\dotsb
s_{2}s_{0}\\
\segmentD{}{1}(j+1)&=s_{t-1}s_{t}s_{t+1}\dotsb
s_{n-2}s_{n-1}s_{n}s_{n-2}\dotsb s_{2}s_{1}.
\end{align*}
By Lemma~\ref{l:commutation.D}, we have the commutation moves
$\segmentD{}{1}(j+1) s_{i} = s_{i} \segmentD{}{1}(j+1) $ for all $t
\leq i\leq n-2$, $ \segmentD{}{1}(j+1) s_{n-1} = s_{n}
\segmentD{}{1}(j+1) $, $\segmentD{}{1}(j+1) s_{n} = s_{n-1}
\segmentD{}{1}(j+1)$.  Note, $\segmentD{}{1}(j+1) s_{t-1}$ is not
equal to a single generator times a segment, but $\segmentD{}{1}(j+1)
s_{t-1} = s_{t-1} \cdot \segmentD{}{1}(j) s_{t-1} = s_{t-1}
s_{t-2}\cdot \segmentD{}{1}(j)$.   Then by Lemma~\ref{l:commutation.D}
again we have $\segmentD{}{1}(j) s_{i} = s_{i-1} \segmentD{}{1}(j) $
for all $2\leq i<t-1$.  Therefore,
\begin{align*}
\segmentD{}{1}(j+1) &\segmentD{}{0}(j) = \segmentD{}{1}(j+1) \cdot
s_{t}s_{t+1}\dotsb s_{n-2}s_{n-1}s_{n}s_{n-2} \dotsb s_{2}s_{0}
\\
= & s_{t}s_{t+1}\dotsb s_{n-2}s_{n-1}s_{n}s_{n-2} \dotsb s_{t} \cdot
\segmentD{}{1}(j+1) \cdot s_{t-1}\dotsb s_{2}s_{0}
\\
= & s_{t}s_{t+1}\dotsb s_{n-2}s_{n-1}s_{n}s_{n-2} \dotsb s_{t} \cdot
s_{t-1} \segmentD{}{1}(j) \cdot s_{t-1}\dotsb s_{2}s_{0}
\\
= & s_{t}s_{t+1}\dotsb s_{n-2}s_{n-1}s_{n}s_{n-2} \dotsb s_{t} \cdot
s_{t-1} \cdot s_{t-2}\dotsb s_{1} \cdot \segmentD{}{1}(j) \cdot s_{0}\\
= & \segmentD{}{1}(j) \cdot \segmentD{}{1}(j) \cdot s_{0}
\\
= & \segmentD{}{1}(j) \cdot s_{t}s_{t+1}\dotsb
s_{n-2}s_{n-1}s_{n}s_{n-2}s_{n}\dotsb s_{2}s_{1} \cdot s_{0}
\\
= &\segmentD{}{1}(j)
\cdot s_{t}s_{t+1}\dotsb s_{n}\dotsb s_{2}s_{0} \cdot s_{1} 
\\
= &\segmentD{}{1}(j) \cdot \segmentD{}{0}(j) \cdot s_{1} 
\end{align*}
since $s_{0}$ and $s_{1}$ commute. 

For $j=n-1$, the proof is similar to the calculations above.
\end{proof}

\begin{theorem}\label{t:surjective.Dn}
In type $D_{n}$, $\pi : \minreps \longrightarrow \partitions =
\partitions(\Dn)$ is a length preserving bijection.
\end{theorem}

\begin{proof} 
As in the proof of Theorem~\ref{t:surjective.Bn}, we only need to show
that the allowed pairs in type $\Dn$ are colored partitions in
$\partitions(\Dn)$ given Lemma~\ref{l:Dpartitions}.  This follows
directly from Lemma~\ref{l:C.commutations.D}.
\end{proof}

\begin{corollary}\label{c:allowed.pairs.D}
The product $\segment{y}{1}(j) \cdot \segment{z}{0}(k) \in \minreps$ and
$\ell(\segment{}{1}(j) \cdot \segment{}{0}(k))=j+k$ if and only if
$y=z$ and ($1\leq j<k\leq 2n-2 $ or $n-1\leq j=k \leq 2n-2$).  
\end{corollary}

\begin{ex}\label{ex:labeling.D}
If $n=6$, the colored partition $(7,5^{b},5^{b},3,1) \in
\partitions(D_{n})$ corresponds with $s_0\cdot s_{3}s_{2}s_{1} \cdot
s_{6}s_{4}s_{3}s_{2}s_{0} \cdot s_{5}s_{4}s_{3}s_{2}s_{1} \cdot
s_{4}s_{6}s_{5}s_{4}s_{3}s_{2}s_{0}$.  Pictorially, we have
\[
\tableau{
0 & 2 & 3 & 4 & 6 & 5 & 4\\
1 & 2 & 3 & 4 & 5 &  & \\
0 & 2 & 3 & 4 & 6 &  & \\
1 & 2 & 3 &  &  &  & \\
0 &  &  &  &  &  & }
\]
\end{ex}
\noindent where again the corresponding reduced expression is read
along the rows from right to left, bottom to top. Note, the 0's and
1's alternate in the first column starting with a 0 in the top left
corner.  The columns indexed by $1<j<n-1$ contain $j$, for $n< j <
2n-2$ column $j$ contains $2n-j-1$, and column $2n-2$ again alternates
from 0 to 1 starting from the top.  Columns $n-1$ and $n$ alternate
between the values $n-1,n$ as well depending on the unique color of
the parts of size $(n-1)$ if they are present.

Bruhat order on $\minreps$ in type $D$ is similar to type $B$ up to rank
$n-1$.  The proofs for the next two lemmas are very similar to the
proofs of Lemmas~\ref{l:youngs.lattice.B} and~\ref{l:chains.B}.

\begin{lemma}\label{l:youngs.lattice.D}
If $\lambda \in \partitions(D_{n})$ and every $\lambda_{k} \leq n-1$,
then the interval in Bruhat order between $\id$ and $\pi ^{-1}
(\lambda)$ in $\minreps$ is isomorphic as posets to the interval from the
empty partition to $\lambda $ in Young's lattice on strict partitions.
\end{lemma}

\begin{lemma}\label{l:chains.D}
If $\lambda \in \partitions(D_{n})$ is a single row, then the interval
in Bruhat order between $\id$ and $\pi ^{-1} (\lambda)$ is isomorphic
as posets to the interval from the empty partition to $\lambda $ in
the generalized Young's lattice on colored partitions where the parts
of size $n-1$ have two possible colors both covering $(n-2)$ and
covered by $(n)$.  In particular,
\[
P_{w}(t) = \sum_{v\leq w}t^{l(w)} = 1 +t+t^{2}+\dotsb +t^{n-2}+2
t^{n-1}+ t^{n}+\dotsb +t^{2n-2}.
\]
\end{lemma}

Consider the elements in $\minreps$ determined by their corresponding
partitions in $\partitions(D_{n})$ for $n\geq 5$.  By
Lemma~\ref{l:youngs.lattice.D}, Bruhat order up to rank 3 on $\minreps$
is the same as Young's lattice on strict partitions of size at most 3
as in \eqref{e:hasse.diag.B}.  However, for $n=4$, there are 3 colored
partitions of size 3 in $\partitions(D_{4})$.

\begin{equation}\label{e:hasse.diag.D}
\begin{array}{c}
\tableau{{0} & {2} &{3}}   \hspace{.2in} \tableau{{0} & {2} &{4}}   \hspace{.2in} \tableau{{0} & {2}\\ {1} & } \\ 
 \backslash \hspace{.2in} |    \hspace{.2in} /   \\
\\
\tableau{{0} & {2}}\\
| \\
\tableau{{0} }\\
| \\
\emptyset
\end{array}
\end{equation}

\begin{theorem}\label{t:rationally.smooth.D} Let $\widetilde{W}$ be
the affine Weyl group of type $D_{n}$ for $n\geq 4$ and let $w \in
\minreps$.  Then $w$ is palindromic if and only if $\pi(w)$ is a one
row shape $(j)$ with $1\leq j\leq n-1$ or $j=2n-2$ or a staircase
shape $(k, k-1, k-2,\dots, 1)$ for $k\leq n-1$.  Note, if $j=n-1$ or
$k=n-1$, there are two distinct elements of either shape corresponding
to the two possible colors of parts of size $n-1$.
\end{theorem}

This proof is very similar to Theorem~\ref{t:rationally.smooth.B}. In
fact it is a bit easier in this case; to show $P_{w}(t)=\sum a_{i}
t^{i}$ is not palindromic it will suffice to show that
$a_{l(w)-j}>1=a_{j}$ for $j=1$ or $j=2$.

\begin{proof}
By Lemma~\ref{l:chains.D}, we know that $P_{w}(t) = 1 +t +\dotsb
+t^{j}$ if $\pi (w)$ is a one row shape with $0\leq j\leq n-1$, which
is clearly palindromic.  Also by Lemma~\ref{l:chains.D}, if $\pi
(w)=(2n-2)$, then $P_{w}(t)$ is palindromic.  By
Lemma~\ref{l:youngs.lattice.D}, if $\pi (w)$ is a staircase shape $(k,
k-1, k-2,\dots, 1)$ for $1<k\leq n-1$ then
\[
P_{w}(t) = \prod_{j=1}^{k} (1+t^{j})
\]
which is palindromic.

Conversely, it follows from Lemma~\ref{l:chains.D} that if $\pi(w)$ is
a one row shape of length $n\leq j<2n-2$, then $P_{w}(t)$ is not
palindromic.  So, assume $\lambda = \pi(w)$ has at least two rows and
is not a staircase shape for $k\leq n-1$.  From
\eqref{e:hasse.diag.D}, we know $a_{1}=a_{2}=1$.  If $\lambda$ has two
or more $\partitions$-removable outside corners, then $a_{l(w)-1} \geq 2$ so $w$
is not palindromic.  Similarly, if $\lambda$ has one $\partitions$-removable
corner and removing it leaves a shape with two $\partitions$-removable corners, then
$a_{l(w)-2} \geq 2$ so $w$ is not palindromic.  One of these
two cases occurs for all partitions $\lambda \in \partitions(D_{n})$
with at least two rows and $\lambda_{1}\geq n$ or if there exists an
$i$ such that $\lambda_{i} -\lambda_{i+1}\geq 2$.  Therefore, it
remains to show that if $\lambda = ((n-1)^{k})$ or $\lambda =
((n-1)^{k},(n-2),\ldots,1)$ for some $k\geq 2$ and with either color
pattern, then $P_{w}(t)$ is not palindromic.  We can fix the color
pattern for the parts of size $n-1$ and look only at partitions with
the exact same color pattern in the proof below.  Therefore, we will
not need to specify the fixed color pattern. 

Say $\lambda = ((n-1)^{k})$ for $2\leq k \leq n-1$.  Then $w$ covers the
two elements $\pi^{-1}((n-1)^{k-1},n-2)$ and $\pi^{-1}(n^{k-1},
n-k-1)$ as seen by omitting the $s_{n-k}$ in the last row of
$\lambda$ and allowing all the larger generators to commute up.  For
example, for $n=6$, we omit the 2 on the last row and allow the
$s_{6}s_{4}s_{3}$ to commute up via the relations in
Lemma~\ref{l:commutation.D}
\[
\tableau{
0 & 2 & 3 & 4 & 5\\
1 & 2 & 3 & 4 & 6\\
0 & 2 & 3 & 4 & 5\\
1 & {} & 3 & 4 & 6
}    
\cong
\tableau{
0 & 2 & 3 & 4 & 5 &6\\
1 & 2 & 3 & 4 & 6 & 5\\
0 & 2 & 3 & 4 & 5 & 6\\
1 &  &  &  &  &
}
\]
Hence, $a_{1}=1 < 2 \leq a_{l(w)-1}$.  

If $\lambda = ((n-1)^{k})$ with $k\geq n$, then $w$ covers
$v=\pi^{-1}((n-1)^{k-1},n-2)$, and $v$ covers two elements since
$n\geq 4$, namely the ones corresponding to the partitions
$((n-1)^{k-1},n-3)$ and $((n+1)^{n-2},(n-1)^{k-n})$.  The ladder
element is obtained by omitting the square in column $n-1$ of row $n$
and the square in column 1 of row $n-1$ of $\lambda$.  As in the case
of the $(n^{n})$ square in the proof of
Theorem~\ref{t:rationally.smooth.B}, the remaining generators in rows
$n$ and $n-1$ commute up to add two additional squares to each of the
first $n-2$ rows.  Therefore, $a_{l(w)-2}\geq 2 > a_{2}=1$.  For
example, when $n=k=4$ we have
\[
\tableau{
0 & 2 & 4\\
1 & 2 & 3\\
{} & 2 & 4\\
1 & 2 & {}\\
}
\cong
\tableau{
0 & 2 & 4 & 3 & 2\\
1 & 2 & 3 & 4 & 2 
 &  & \\
 &  & \\
}
\]

On the other hand, say $\lambda = ((n-1)^{k},(n-2),\ldots,1)$ for some
$k\geq 2$.  If $k\geq n$, then we can delete two generators in
$\lambda$ in two ways as in the proof of
Theorem~\ref{t:rationally.smooth.B} for the shape $((n-1)^{n})$ and
obtain two elements below $w$ of length $l(w)-2$.  If $k=2$, then
there are two elements of length $l(w)-2$ below $w$ corresponding to
the partitions $((n-1)^{2},(n-2),\dotsc,3,1)$ and $(n,n-1,n-2,\dotsc
,3)$ obtained from $\lambda$ by deleting the last two elements in
column 1 and allowing the $s_{2}$ to commute up.  If $2< k<n$, then
again $w$ covers two elements corresponding to the partitions
$((n-1)^{k},(n-2),\ldots,3,2)$ and $(n^{k-2},(n-1)^{2},n-2,\dotsc,
\widehat{k-1},\dotsc,1)$, the ladder element is obtained from
$\lambda$ as follows.  Delete the largest generator $s_{k-1}$ in row
$n$.  The corresponding product of segments has two consecutive
segments of length $k-2$, use the segment relations
$\segment{}{1}(k-2)\segment{}{0}(k-2)=
\segment{}{1}(k-3)\segment{}{0}(k-2)\cdot s_{1}$ or the analogous
relation switching $0$ and $1$.  Then the $s_{1}$ commutes up to row 1
and becomes an $s_{n}$ or $s_{n-1}$ extending row 1.  After applying
the segment relation, now rows $n+1$ and $n+2$ have the same length.
Repeating the same procedure adds an extra box to row 2, etc.  After
applying $k-2$ segment relations, the remaining partition is
$(n^{k-2},(n-1)^{2},n-2,\dotsc, \widehat{k-1},\dotsc,1)$.
\end{proof}

\begin{corollary}\label{c:mw.D}
For $W$ of type $D_{n}$, we have $m_{W}=n-1$.  
\end{corollary}

\begin{proof}
The proof above shows that $n-1$ is sufficient.  It is necessary for
the element $w=\segment{}{0}(2n-3) \in \minreps$ with $P_{w}(t)= t^{n}+\sum_{i=0}^{2n-3}t^{i}$.
\end{proof}

\section{Exceptional Types}\label{s:exceptional}

In this section we give the algorithm we used to prove
Theorem~\ref{t:partition.bijection} for the exceptional types.  As an example of the algorithm,
we give the full proof in type $G_{2}$.  In type $F_{4}$, we give the
necessary data to verify the theorem without much discussion in the
appendix.  For
types $E_{6},E_{7},E_{8}$ we refer the reader to our Lisp and Maple code available at 
\begin{center}
\textrm{http://www.math.washington.edu/$\sim$billey}
\end{center}

Every exceptional type has a Type I Coxeter graph.  The segments can
be obtained using the Mozes numbers game
\cite{b-b,eriksson}.  The node-firing game from
\cite{Billey-Mitchell} can be used as well since the segments form an
interval in left weak order so one only needs to find the longest
element in $\minreps$ with exactly one $s_{0}$ in any reduced
expression.  The allowed pairs can also be identified using the
numbers game or the node-firing game.  Below is the high level
algorithm for verifying Theorem~\ref{t:partition.bijection}:

\begin{enumerate}
\item Identify the set of all segments using the
node-firing game.
\item Compute all products of pairs $\segment{c}{}(i)
\segment{d}{}(j)$ to identify allowed pairs: \\
$\segment{c}{}(i)
\segment{d}{}(j)$ is \textit{allowed} if and only if
\[
\ell(\segment{c}{}(i) \segment{d}{}(j)) = i+j
\]
and 
\[
\segment{c}{}(i) \segment{d}{}(j) \in \minreps.
\]
Both tests are straightforward using the numbers game.  

\item Let $\mathcal{A}$ be the finite partial order on the allowed
colored parts in $\partitions$ given by $i^{c} < j^{d}$ if $(i^{c},
j^{d})$ is an allowed pair.  Verify that $\mathcal{A}$ is a subposet
of left weak order restricted to the segments confirming
Lemma~\ref{l:pairs}.\ref{l:weak.order.and.segmetns} in this case.

\item Construct a representation of the rational generating function
$G_{\partitions} = \sum_{\lambda \in \partitions} t^{|\lambda|}$
corresponding with the affine partitions.  If $(j^{d},j^{d})$ is an
allowed pair, we say $j^{d}$ is \textit{repeatable}.  Let
\[
D(j^{c}) =
\begin{cases}
(1 - t^{j}) &	j^{c} \text{ repeatable}\\
1 & \text{otherwise }
\end{cases}
\]
and if $\lambda = (\lambda_{1}^{c_{1}}>\lambda_{2}^{c_{2}}>\dotsb >
\lambda_{k}^{c_{k}})$ is a strict colored partition, let $D(\lambda ) =
\prod_{i=1}^{k}D(\lambda_{i}^{c_{i}})$.  Then

\begin{equation}\label{e:p.expansion}
G_{\partitions} = \sum \frac{t^{|\lambda|} }{D(\lambda)}
\end{equation}
where the sum is over all strictly decreasing chains $\lambda =
(\lambda_{1}^{c_{1}}>\lambda_{2}^{c_{2}}>\dotsb > \lambda_{k}^{c_{k}})$
(including the empty chain) in $\mathcal{A}$ and $|\lambda|= \sum
\lambda_{i}$.  The sum in \eqref{e:p.expansion} is finite since there
are only a finite number of segments in each type.  

\item Pipe the rational generating function $G_{\partitions}$ into
Maple and use the command ``simplify($G_{\partitions}$)''.  Compare
the result with Bott's formula.
\end{enumerate}

\begin{rem}
In practice a reasonably efficient way to compute the rational
generating function in \eqref{e:p.expansion} is to identify the set of
all strict colored partitions with no repeatable parts and then for
each subset of the repeatable parts, check if the repeatable parts can
appear with the non-repeatable parts.  Furthermore, the strictly
decreasing chains of non-repeatable parts should be constructed using a
depth first search that only looks down branches corresponding with
allowed pairs.
\end{rem}

\noindent \textbf{Example in Type $G_{2}$}: The Dynkin diagram of type
$\widetilde{G}_{2}$ is a Type I graph.  The segments are as follows:

\begin{equation}\label{e:segments.g2}
\begin{array}{lc}
\text{segment } & \text{length }
\\
\hrulefill & \hrulefill 
\\
\segment{}{}(1)= s_0&  1
\\
\segment{}{}(2)= s_1 s_0&  2
\\
\segment{}{}(3)= s_2 s_1 s_0&  3
\\
\segment{}{}(4)=s_1 s_2 s_1 s_0&  4
\\
\segment{}{}(5)=s_2 s_1 s_2 s_1 s_0&  5
\\
\segment{}{}(6)=s_1 s_2 s_1 s_2 s_1 s_0&  6
\end{array}
\end{equation}

Observe that the segments have distinct lengths.  No colors will be
necessary for the corresponding partitions.

One can verify that the product of two segments $\segment{}{}(i)\cdot
\segment{}{}(j)$ is a reduced expression for a minimal length coset
representative in $\minreps$ if and only if the ordered pair $(i, j)$
appears in the following list of \textit{allowed pairs}:
\[
\{(1, 4), (1, 5), (1, 6), (2, 5), (2, 6), (3, 5), (3, 6), (4, 5), (4, 6), (5, 5), (5, 6), (6, 6) \}
\]
Let $\partitions(G_{2})$ be the set of partitions allowed by these
pairwise rules, i.e. they have an arbitrary number of 5's and 6's
followed by one of the partitions from the following set:
\[
\{(1),(2),(3),(4),(4,1) \}
\]
The generating function for such partitions is 
\begin{equation}\label{e:generating.function.g2}
\frac{1+t+t^{2}+t^{3}+t^{4}+t^{1}t^{4}}{(1-t^{5})(1-t^{6})} =
\frac{1+t+t^{2}+t^{3}+t^{4}+t^{5}}{(1-t^{5})(1-t^{6})} = \frac{1}{(1-t^{1})(1-t^{5})}
\end{equation}
Observe that the formula above agrees with Bott's formula for the rank
generating function for $\minreps$ since the exponents for type
$G_{2}$ are $\{1,5 \}$.  This proves
Theorem~\ref{t:partition.bijection} for type $G_{2}$.

Above it was easy enough to use the description of
$\partitions(G_{2})$ to get the generating function in
\eqref{e:generating.function.g2} since 5 and 6 can be added to every
affine partition to get another affine partition.  In types $F_{4},
E_{6,7,8}$ some repeatable parts are not allowed to appear with other
parts so we need to be more careful.  For example, in type $F_{4}$ the
part denoted $11^{1}$ in the appendix cannot appear with all other
parts.  Therefore, it is instructive to build this generating function
using \eqref{e:p.expansion} as an example.  The following are all
strict affine partitions in $\partitions(G_{2})$ and their
corresponding denominators:
\[
\begin{array}{lll}
\{(), (1),(2), (3), (4), (1,4) \}  &&	 D(\lambda) = 1
\\
\{(5) \} \cup \{(j,5) \given 1\leq j\leq 4 \}  \cup \{(1,4,5) \}  &&	 D(\lambda) = (1-t^{5})
\\
\{6 \}  \cup\{(j,6) \given 1\leq j\leq 4 \} \cup \{(1,4,6) \}  &&	 D(\lambda) = (1-t^{6})
\\
\{(5,6) \} \cup \{(j,5,6) \given 1\leq j\leq 4 \} \cup \{(1,4,5,6) \}  &&	 D(\lambda) = (1-t^{5})(1-t^{6})
\end{array}
\]
Putting this information into \eqref{e:p.expansion} we have
\[
G_{\partitions} = \left(1+t+t^{2}+t^{3}+t^{4} +t^{1}t^{4} \right)
\left( 1+ 
\frac{t^{5}}{(1-t^{5})} +\frac{t^{6}}{(1-t^{6})}
+\frac{t^{5}t^{6}}{(1-t^{5})(1-t^{6})} \right)
\]
which simplifies to \eqref{e:generating.function.g2} as expected.

\section{Type $A$}\label{s:type.a}

Fix $W$ to be the Weyl group of type $A_{n}$ and $\widetilde{W}$ its
affine Weyl group.  The elements in $\minreps$ are known to be in
bijection with $k+1$-core partitions and $k$-bounded partitions
\cite{EE-98,Misra-Miwa,Lapointe-Morse-2005}.  In this section, we will
describe the affine partitions in type $A$ and prove that these
objects are also in bijection with $\minreps$.    The proof relies on
the $k+1$-core bijection.  We use the segments to identify the
palindromic elements as well.  Following a suggestion due to Hugh
Thomas, we use the $k$-bounded partitions to identify the Poincar\'e
polynomial for all the palindromic, nonsmooth elements in $\minreps$
reproving a theorem due to the second author \cite{Mitchell-86}.

In the Dynkin diagram for $\widetilde{A}_{n}$, both $s_{1}$ and
$s_{n}$ are connected to $s_{0}$ in the Dynkin diagram for
$\widetilde{W}$.  We will treat $W$ as a Type I Coxeter group with
$J=\langle s_{2},\dots, s_{n-1} \rangle $.  The fragments in this type are the
elements of $W^{J}$.  As permutations, these fragments have one line
notation in $S_{n+1}$ of the form
\[
[a, 1,2,\dots, \hat{a},...,\hat{b},\dots, n+1,b] = s_{a-1}s_{a-2}\dotsb s_{1}s_{b}s_{b+1}\dotsb s_{n-1}s_{n}
\]
where $a,b \in \{1,\dots, n+1 \}:=[n+1]$ and $a \neq b$.  The segments can be
defined by
\[
\segment{c}{}(k) = s_{c}s_{c-1}\dots s_{1} 
s_{n-k+c+2} \dots s_{n-1}s_n s_{0}
\]
for $0\leq c\leq n$, $1\leq k-c\leq n$, so $1\leq k \leq 2n$ and $k$
is the length of $\segment{c}{}(k)$.  There are $n(n+1)$ segments for
type $A_{n}$.

To describe the weak order on segments and the allowed pairs, it is
convenient to use an alternative notation.  For $0\leq i\leq n$ and
$1\leq j\leq n$, let $C_{i,j} = s_{i}s_{i-1}\dots s_{1} s_{n-j+2}
\dots s_{n-1}s_n s_{0}$ so $\ell(C_{i,j})=i+j$ and $\segment{c}{}(k) =
C_{c,k-c}$.

\begin{lemma}\label{l:relations.An}
For $1\leq a\leq n$,  we have the following relations:
\[
s_{a}C_{i,j} = \begin{cases}
C_{i,j}s_{a+1} &	 a<i \text{ and } a < n-j \\
C_{i,j+1} &        	 a<i \text{ and } a = n-j \\
C_{i,j-1} &        	a<i \text{ and } a = n-j+1 \\
C_{i,j}s_{a} &        	 a<i \text{ and } a \geq   n-j+2 \\
C_{i-1,j}   &	        a=i\\
C_{i+1,j}   &	        a=i+1\\
C_{i,j}s_{a} &        	a>i+1 \text{ and } a \leq   n-j \\
C_{i,j+1}  &        	a>i+1 \text{ and } a =   n-j+1 \\
C_{i,j-1}  &        	a>i+1 \text{ and } a =   n-j+2 \\
C_{i,j}s_{a-1} &       	a>i+1 \text{ and } a >  n-j+2 \\
\end{cases}
\]
\end{lemma}

\begin{proof}
These relations follow directly from the Coxeter relations.  
\end{proof}

As a direct corollary of the relations above, we get the following
characterization of the weak order on segments.

\begin{lemma}\label{l:left.order.a}
In left weak order, 
\[
C_{i,j}\leq C_{k,l} \iff  
\begin{cases}
i\leq k \text{ and } j\leq l   &	\text{for $i + j \neq n.$ }
\\
i<k\text{ and } j\leq l    &	\text{for $i+ j=n.$ }
\end{cases}
\]
\end{lemma} 
Hence, the left weak order on segments is a poset similar to the
product of chains $[n+1] \times [n]$, except that it is missing some
relations between the middle two ranks namely, $C_{i,j} \not <
C_{i,j+1}$ if $i+j=n$.

Thus, the weak order on segments in type $A$ is significantly more
complicated than in all other types.  In order to characterize the
allowed pairs $C_{ij}C_{k,l}$ and finish the proof of
Theorem~\ref{t:partition.bijection} in type $A$, we introduce the
bijection from $\minreps$ to core partitions.   

\subsection{Core Partitions, $k$-bounded
partitions, and affine partitions}

Let $\lambda$ be a partition thought of as a Ferrers diagram.  We say
a square $s$ is \textit{addable} to $\lambda$ if $s$ is adjacent to
the southeast boundary of $\lambda$ and adding $s$ to $\lambda$
results in a larger partition shape.  Similarly, a square $s$ inside
$\lambda$ is \textit{removable} from $\lambda$ if $s$ is adjacent to
the southeast boundary of $\lambda$ and removing $s$ from $\lambda$
results in a smaller partition shape.  The \textit{$n$-content} of the
square $(i,j)$ in matrix notation is $(j-i) \mathrm{mod}(n+1)$.  The
\textit{arm} of a square $s$ in $\lambda$ is the set of boxes strictly
to the right of $s$.  The \textit{leg} of a square $s$ in $\lambda$ is
the set of boxes strictly below $s$.  The \textit{hook length} of
square $s=(i,j)$ in $\lambda$, denoted $h(s)$, is the number of
squares weakly below $s$ plus the number of squares strictly to the
right of $s$, so $h(s)=|arm(s)|+|leg(s)|+1$.  A \textit{central hook}
of a partition is any hook whose corner square $s$ is along the main
diagonal so $s=(a,a)$ for some $a \in \mathbb{N}$.  The \textit{lowest
central hook} of a partition is the central hook with $a$ as large as
possible.

An $n+1$-\textit{core} is a partition with no hook length divisible by
$n+1$, or equivalently, no hook lengths exactly $n+1$
\cite{M1}. In type $A_{n}$, Misra and Miwa \cite{Misra-Miwa}
showed that $n+1$-core partitions are in bijection with $\minreps$.
Lascoux showed that Bruhat order on $\minreps$ in type $A_{n}$ is
completely determined by the partial order of containment on the
corresponding $n+1$-cores \cite{Lascoux-99}.  See also
\cite{bjorner-brenti-1996}.

Let $C^{n+1}$ be the set of all $n+1$-cores. The Misra-Miwa bijection
from
\[
\core: \minreps \longrightarrow C^{n+1}
\]
is defined recursively as follows.  Consider all squares $(i,j) \in
 \mathbb{N} \times \mathbb{N}$ labeled by their $n$-content.  The empty
partition, denoted $\emptyset$, represents the identity element in $\minreps$.  Say
$\core(w)$ is the $n+1$-core for $w\in \minreps$ and $s_{i}w \in
\minreps$ with $l(w)<l(s_{i}w)$.  To obtain the Ferrers diagram for
$\core(s_{i}w)$ from the Ferrers diagram for $\core(w)$, add all of
the addable squares of $\core(w)$ with $n$-content $i$ to $\core(w)$.

For example, if $n =3$, then the $4$-core for $s_{2}s_{3}s_{0}$ is a
single column with three boxes obtained by first adding the unique box
with $n$-content $0$, then $3$, then $2$ in such a way as to obtain a
partition shape at each step.  The $n+1$-cores starting with $s_{2}s_{3}s_{0}$ and
building up to $s_{0} s_{3} s_{2}s_{1}s_{2}s_{3}s_{0}$ are given by
the following shapes (the $n$-content is written inside each box just for
convenience):
\[
\begin{array}{cccccc}
s_{2}s_{3}s_{0} & s_{1}s_{2}s_{3}s_{0} & s_{2}s_{1}s_{2}s_{3}s_{0} &
s_{3}s_{2}s_{1}s_{2}s_{3}s_{0} & s_{0} s_{3}s_{2}s_{1}s_{2}s_{3}s_{0}
&\\
\tableau{\\
{0}   \\
{3}  \\
{2}  \\
  \\
  \\
}
&
\tableau{\\
{0} &{1} &  \\
{3}  & & \\
{2}  & &  \\
{1} & &   \\
&&   \\
}
&
\tableau{\\
{0} &{1} &{2}  \\
{3}  & & \\
{2}  & &  \\
{1} & &   \\
&&   \\
}
&
\tableau{\\
{0} &{1} &{2} & {3}&  \\
{3}  & & \\
{2}  & & & \\
{1}  & & & \\
  & & & \\
}
&
\tableau{\\
{0} &{1} &{2} & {3}& {0} \\
{3} & {0}& & & & \\
{2} & & & & & & \\
{1} & & & & & & \\
{0} & & & & & & \\
}
\end{array}
\]
Note, that the number of boxes in $\core(w)$ is at least $\ell(w)$ but
not necessarily equal to $\ell(w)$.

More formally, for $0\leq i\leq n$, let $r_{i}$ be the operator on
partitions which acts on $\lambda$ by adding all addable boxes of
$n$-content $i$ and removing all removable boxes of $n$-content $i$. If $w\in
\minreps$ and $w=s_{i_{1}}s_{i_{2}}\dots s_{i_{p}}$ is a reduced
expression, then define $\core(w) = r_{i_{1}}r_{i_{2}}\dots
r_{i_{p}}(\emptyset)$.  To see that $\core: \minreps \longrightarrow
C^{n+1}$ is a proper map and a bijection we rely on the following
facts which are proved in \cite{Lapointe-Morse-2005}.

\begin{theorem}\label{t:core.facts} Let $w \in \minreps$.   
\begin{enumerate}
\item\label{i:core.facts.1} If $\lambda \in C^{n+1}$ then $r_{i}(\lambda ) \in C^{n+1}$ for
all $0\leq i \leq n$ so $\core(w) \in C^{n+1}$. 

\item \label{i:core.facts.2} We have $s_{i}w \in \minreps$ if and only
if $\core(w)$ has at least one addable or one removable square of
$n$-content $i$.  If $w <s_{i}w$, then $\core(w)$ has no removable squares
of $n$-content $i$ and all the addable squares of $n$-content $i$ lie in
consecutive diagonals. If $w >s_{i}w$, then $\core(w)$ has no addable
squares of $n$-content $i$ and all the removable squares of $n$-content $i$
lie in consecutive diagonals.  If $w \in \minreps$ and $s_{i}w \not
\in \minreps$, then $\core(w)$ has no addable or removable squares
with $n$-content $i$.

\item \label{i:core.facts.3} The relations among the set of generators
$\widetilde{S}=\{s_{0},s_{1},\dotsc, s_{n} \}$ acting on $\minreps$
from the left are exactly the same as the relation among the operators
$\{r_{0},\dots, r_{n} \}$ acting on $C^{n+1}$.  Thus, we can consider
the action $s_{i}:C^{n+1} \longrightarrow C^{n+1}$ to be same as
$r_{i}$.
\end{enumerate}
\end{theorem}


Now we can return to products of segments.  Roughly speaking, pairwise
allowed products of segments in type $A$ correspond with consecutive central
hooks in the core partition.

\begin{lemma}\label{l:allowed.a}
The product $C_{i,j}C_{k,l}$ is reduced and in $\minreps$ if and only if at least one of the following hold:
\begin{enumerate}
\item $i<k$ and $j<l$. 
\item $k+l>n$ and $i<k$ and $j\leq l$. 
\item $i+j>n$ and $i\leq k$ and $j\leq l$. 
\end{enumerate}
\end{lemma}

\begin{proof}

By Theorem~\ref{t:core.facts}(~\ref{i:core.facts.2}), if $w \in
\minreps$ and $s_{i}\core(w)$ adds at least one box, then $s_{i}w \in
\minreps$, $\ell(s_{i}w) = \ell(w)+1$, and $\core(w) \subset
s_{i}\core(w)$.  In particular, this shows that every segment
$C_{k,l}=s_{k}\dotsb s_{1}s_{n-l+2}\dotsb s_{n}s_{0}$ has the
properties: $C_{k,l}\in \minreps$, $\ell(C_{k,l})=k+l$ and
$\core(C_{k,l})$ is a single hook shape core partition $\core(C_{k,l})$
with hook length at least $k+l$ and arm length exactly k. The
generators $s_{n-l+2}\dotsb s_{n}s_{0}$ build up the leg of the hook,
and the generators $s_{k}\dotsb s_{1}$ build up the arm of the hook.
If $k+l>n$ then $\core(C_{k,l})$ has hook length $k+l+1$ since applying
the generators $s_{k}\dotsb s_{1}$ along the arm of the hook must
necessarily also add one more box to the leg of the hook.

Let $H_{0}=\core(C_{k,l}), H_{1}=s_{0} \core(C_{k,l}), \dotsc ,
H_{i+j} = s_{i}\dotsb s_{1}s_{n-j+2}\dotsb s_{n}s_{0} \core(C_{k,l})$
be the result of applying each of the generators of $C_{ij}=
s_{i}\dotsb s_{1}s_{n-j+2}\dotsb s_{n}s_{0}$ from right to left on
$\core(C_{k,l})$.  Note, $H_{1}$ is a hook plus one extra square on the
central diagonal and $H_{i+j}$ is the union of exactly two central
hooks since there are no other $s_{0}$ generators in the reduced
expression for $C_{i,j}$.

Again by Theorem~\ref{t:core.facts}(~\ref{i:core.facts.2}), the
product $C_{i,j}C_{k,l}$ is reduced and in $\minreps$ if and only if
$H_{a}$ contains at least one square not in $H_{a-1}$ for each $1\leq
a\leq i+j$.  Observe, that if one of the three conditions above are
satisfied, the latter condition holds since the second central hook is
built up exactly as it is when applying the generators of $C_{ij}$ to
the empty partition.

Conversely, assume that none of the three conditions hold on
$i,j,k,l$. Then there exists a smallest $a$ such that either the
lowest squares in both the first and second central hooks of
$H_{a-1}$ are in the same row and the next generator tries to add one
more square to the leg of the second hook but fails, or the rightmost
squares in both the first and second central hooks of $H_{a-1}$ are
in the same column and the next generator tries to add one more square
to the arm of the second hook but fails.  In either case, we claim no
new squares are added in $H_{a}$ so $H_{a} \subset H_{a-1}$ hence
$C_{i,j}C_{k,l}$ cannot be a reduced factorization for an element in
$\minreps$.  

To prove the claim $H_{a} \subset H_{a-1}$, assume the lowest squares
in both the first and second central hooks of $H_{a-1}$ are in the
same row and the next generator tries to add one more square in row
$r$ and column $2$ to the leg of the second hook.  The other case is
similar. The generator being applied to go from $H_{a-1}$ to $H_{a}$
is $s_{n+1+2-r}$ where $n+1+2-r\geq n-j+2>1$ by definition of the
segment $C_{ij}$.  Let $z=n+1+2-r$ so the $n$-content of $(r,2)$ is $z \
\mathrm{mod} (n+1)$.  If the cell $(r,2)$ is not addable then it must
be because $(r,1)$ is not a square in $H_{a-1}$, therefore no square
in any row below row $r$ is addable.  Square $(2,3)$ has $n$-content 1 so
it is not is not added from $H_{a-1}$ to $H_{a}$ since $z>1$, nor is
any square south or east of $(2,3)$.  Finally, if a cell of $n$-content
$z$ was addable in the first row, then the hook with corner $(1,2)$
would have hook length $n+1$ in $H_{a}$, but this cannot happen by
Theorem~\ref{t:core.facts}(~\ref{i:core.facts.1}).  Therefore, $H_{a}$
is contained in $H_{a-1}$.
\end{proof}

Rephrasing the vocabulary slightly from Section~\ref{s:canonical}, we
say two segments $C_{i,j}C_{k,l}$ form an \textit{allowed pair} if the
conditions of Lemma 3 hold.  Note, in the language of
Section~\ref{s:canonical}, we would have said $((i+j)^{i}, (k+l)^{k})$
is an allowed pair as a colored partition since
$C_{ij}=\segment{i}{}(i+j)$.  However, in type $A$ it is more
convenient to work with the $C_{ij}$ notation for segments.
Furthermore, we say $C_{i_{1},j_{1}}C_{i_{2},j_{2}}\dotsb
C_{i_{p},j_{p}}$ is an \textit{allowed product of segments} if each
consecutive pair is an allowed pair.

\begin{theorem}\label{t:surjective.An}
Let $W$ be the Weyl group of type $A_{n}$.  There is a length
preserving bijection from $\minreps$ to allowed product of segments
given by the canonical factorization $r(w)=
C_{i_{0},j_{0}}C_{i_{1},j_{1}}\dotsb C_{i_{p},j_{p}}$.  
\end{theorem}

\begin{proof}
Injectivity of the map $r$ is proved in the first paragraph of the
proof for Theorem~\ref{t:partition.bijection}.  Therefore, it remains
to show that any allowed product of segments gives rise to a reduced
expression for an element in $\minreps$.

Let $v = C_{i_{1},j_{1}}\dotsb C_{i_{p},j_{p}}$ be an allowed product
of segments. By induction we can assume $C_{i_{1},j_{1}}\dotsb
C_{i_{p},j_{p}}$ is a reduced factorization for $v$ and $v \in
\minreps$.  Assume $C_{i_{0},j_{0}}C_{i_{1},j_{1}}$ is also an allowed
product.  We will show that $w=C_{i_{0},j_{0}}v \in \minreps$ and
$l(w)=l(v)+i_{0}+j_{0}$.

Construct the $n+1$-core partition $\core(v)$ by multiplying one
segment at a time, starting with $C_{i_{p},j_{p}}$ and ending with
$C_{i_{1},j_{1}}$.  By induction, we assume that multiplication by
$C_{i_{1},j_{1}}$ in this construction adds a new square on the main
diagonal with $n$-content 0 in square $s=(a,a)$ for some $a$, new squares
with $n$-content $n,n-1,\dotsc, n-j_{1}+2$, below $s$, and then new
squares with $n$-content $1,2,\dotsc,i_{1}$ to the right of $s$, plus
possibly other squares, thus the length of the corresponding element
in $\minreps$ is strictly increasing with each generator.  Similarly,
since $C_{i_{0},j_{0}}C_{i_{1},j_{1}}$ is also an allowed product,
starting with $\core(v)$ and multiplying by the generators in
$C_{i_{0},j_{0}}$ in reverse order will add at least one new box each
with each generator as shown in Lemma~\ref{l:allowed.a}.  The
resulting partition will be an $n+1$-core $\lambda$.  Therefore, the
expression $C_{i_{0},j_{0}}C_{i_{1},j_{1}}\dotsb C_{i_{p},j_{p}}$ is a
reduced factorization for $w = \core^{-1}(\lambda)$ so $w \in
\minreps$.
\end{proof}

\begin{corollary}\label{c:p.central.hooks}
Assume $w \in \minreps$.  The factorization $r(w) =
C_{i_{1},j_{1}}C_{i_{2},j_{2}}\dotsb C_{i_{p},j_{p}}$ can be obtained
from $\core(w)$ by successively removing the lowest central hook of
the remaining partition for each factor.  To remove a hook, apply
generators which remove squares in the lowest central hook along the
arm first, then the leg, then the corner square.
\end{corollary}

\begin{corollary}\label{c:p.distinct.hooks}
If $r(w) = C_{i_{1},j_{1}}C_{i_{2},j_{2}}\dotsb C_{i_{p},j_{p}}$ then
$\core(w)$ has $p$ squares along the main diagonal.  Conversely, if
$\core(w)$ has exactly $p$ squares along the main diagonal, then there
are $p$ segments in the canonical factorization $r(w)$.
\end{corollary}

\begin{corollary}\label{c:involution}
Let $i:\minreps \longrightarrow \minreps$ be the involution determined
by mapping $s_{i}$ to $s_{n+1-i}$ for $1\leq i \leq n$.  If $w=
C_{i_{1},j_{1}}\dotsb C_{i_{p},j_{p}}$, then $i(w) =
C_{j_{1}-1,i_{1}+1}\dotsb C_{j_{p}-1,j_{p}+1}$.
\end{corollary}

\begin{proof}
The proof follows directly from the observation that $\core(w)$ is the transpose of $\core(i(w))$. 
\end{proof}




Let $B^{n}$ be the set of $n$-bounded partitions, so $\lambda \in
B^{n}$ implies $\lambda_{i}\leq n$ for all $i$.  For a fixed $n$,
define a map 
\[
\bounded : C^{n+1} \longrightarrow B^{n}
\]
given by mapping $\lambda \in C^{n+1}$ to $\mu \in B^{n}$ if the
Ferrers diagram for $\lambda$ has $\mu_{i}$ squares with hook length
at most $n$ in row $i$ for each $i\geq 1$.  For example, if $n=3$,
$\lambda =(5,2,1,1,1)$ then $\bounded(\lambda)=(3,1,1,1,1)$.  

Define a second map, 
\[
\boundedtominreps : B^{n} \longrightarrow \minreps
\]
as follows.  Given $\mu \in B^{n}$, fill the squares in the Ferrers
diagram of $\mu$ by their $n$-content.  Say $i_{1}i_{2}\dotsc i_{p}$
is the word obtained by reading along the rows of $\mu $ from
right to left, bottom to top.  Then
$\boundedtominreps(\mu)=s_{i_{1}}s_{i_{2}}\cdots s_{i_{p}}$.  For
example, with $n=3$, 
\[
(3,3) = \tableau{
0 & 1 & 2\\
3 & 0 & 1
}
\hspace{.2in}
\implies
\hspace{.2in}
\boundedtominreps(3,3) =
s_{1}s_{0}s_{3}s_{2}s_{1}s_{0}.
\]

The next theorem follows from \cite[Lemma 23]{EE-98}, \cite[Theorem
7]{Lapointe-Morse-2005}, and \cite{Misra-Miwa}.

\begin{theorem}\label{t:lapoint-morse}
The maps $\boundedtominreps$ and $\bounded$ are bijections such that
$\boundedtominreps \circ \bounded \circ \core : \minreps
\longrightarrow \minreps$ is the identity map.  Furthermore, if
$\boundedtominreps(\mu ) = w \in \minreps$ then $|\mu| = \ell(w)$.
\end{theorem}

\begin{corollary}\label{c:affine.part.to.bounded}
In type $A_{n}$, there exists a size preserving bijection from
$n$-bounded partitions to affine partitions $\partitions$.
\end{corollary}

Note, in general the canonical reduced expressions given by the
$n$-bounded partition and the corresponding affine partition are
different. For example, with $n=3$ again

\[
(2,2,1) = \tableau{
0 & 1 \\
3 & 0 \\
2 &
}
\hspace{.2in}
\implies
\hspace{.2in}
\boundedtominreps(2,2,1) =
s_{2}s_{0}s_{3}s_{1}s_{0} = s_{0}s_{2}s_{1}s_{3}s_{0} = C_{0,1}C_{2,2}.
\]

\subsection{Palindromy in type $A$}\label{s:pal.a}

Given that the partition bijection $\pi : \minreps \longrightarrow
\partitions$ holds in type $A$ by Theorem~\ref{t:partition.bijection}
and Theorem~\ref{t:surjective.An}, we can identify all the thin
elements of $\minreps$ in type $A$ and cpo's relatively easily.
Furthermore, we can use the relations in Lemma~\ref{l:relations.An} to
identify which of the thin elements are not palindromic.  Therefore,
we can give a complete characterization of palindromic elements in
type $A_{n}$ using the following theorem due to the second author.
This allows us to determine $m_{W}$ as well.

Note, that in type $A_{1}$, $\minreps$ is a chain under Bruhat order
so every element is thin, extra thin and palindromic. Therefore,
throughout the rest of this section we will assume $n\geq 2$.

\begin{theorem}\label{t:mitchell}\cite{Mitchell-86}
In type $A_{n}$ for $n\geq 2$, there exists two infinite families of
palindromic elements, namely $\{C_{0,j}(C_{1,n})^{k} : j+k = 0\
\mathrm{mod}\ n \}$ and $\{C_{i,1}(C_{n,1})^{k} : i+k+1 = 0 \ 
\mathrm{mod}\ n \}$.   
\end{theorem}

We call the palindromic elements in Theorem~\ref{t:mitchell}
\textit{spiral elements} since their unique reduced expressions spiral
around the Dynkin diagram clockwise or counterclockwise.

\begin{lemma}\label{l:thin.a}
Assume $W$ has type $A_{n}$ for $n\geq 2$.  

\begin{enumerate}
\item For each segment $C_{i,j}$ such that $i+j\leq n$ and each $k\geq
0$, the following elements are thin:
\begin{equation*}
\begin{array}{ll}
C_{0,j-i}C_{1,j-i+1}\dotsb C_{i-1,j-1}\cdot C_{i,j} \cdot  (C_{i+1,n-j})^{k} &  \text{if  }0\leq  i<j\leq n, \\
C_{0,j-i}C_{1,j-i+1}\dotsb C_{i-1,j-1}\cdot C_{i,j} \cdot (C_{n+1-i,j})^{k}&  \text{if  } 0\leq  i<j\leq n,\\
C_{i-j,1} C_{i-j+1,2} \dotsb C_{i-1,j-1}\cdot C_{i,j} \cdot (C_{i+1,n-j})^{k} &  \text{if  } 1\leq j\leq i < n, \\
C_{i-j,1} C_{i-j+1,2} \dotsb C_{i-1,j-1} \cdot C_{i,j} \cdot (C_{n+1-i,j})^{k} & \text{if  } 1\leq j\leq  i < n.
\end{array}
\end{equation*}
\item For each pair $1\leq i,j\leq n$ such that $i+j=n+1$ and $k\geq 0$, the elements 
$(C_{i,j})^{k} $ are thin. 
\end{enumerate}
Conversely, every thin element in $\minreps$ is listed above.
\end{lemma}

\begin{proof}
Observe that the branching number $b_{W}=2$, so an affine partition
$\lambda = \pi (w) \in \partitions$ is thin if and only if $\lambda$
covers one element in the generalized Young's lattice.

The only thin segments are of the form $C_{i,j}$ for $i+j=n+1$,
$C_{0,j}$ for $1\leq j\leq n$, and $C_{i,1}$ for $0\leq i< n$.  All
other segments cover two elements given by reducing the first index or
the second index.  Furthermore, the segments $C_{i,j}$ such that
$i+j>n+1$ all cover two segments $C_{i-1,j}$ and $C_{i,j-1}$ and both
are repeatable.  So, if $w$ has any segment of length strictly longer
than $n+1$ then $w$ is not thin.

Therefore, $w=C_{i_{1},j_{1}} \cdots C_{i_{p},j_{p}} \neq \mathrm{id}$
is thin if and only if all of the following hold:
\begin{enumerate}
\item $i_{1}=0$ or $j_{1}=1$ or $i_{1}+j_{1}=n+1$.
\item For each $1\leq k\leq p$, we have $i_{k} + j_{k} \leq n+1$.
\item For each $1\leq k<p$, the pair $(i_{k+1},j_{k+1})$ is minimal so
that $C_{i_{k},j_{k}}C_{i_{k+1},j_{k+1}}$ is an allowed pair.
\end{enumerate}
The lemma now follows from these three conditions.
\end{proof}

\begin{lemma}\label{l:cpos.a}
For each segment $C_{i,j}$ such that $i+j\leq n$ there exists a cpo
with $C_{i,j}$ as its longest segment given by the product

\begin{equation}\label{e:cpo.type}
\begin{array}{ll}
C_{0,j-i}C_{1,j-i+1}\dotsb C_{i-1,j-1}C_{i,j} & \text{ if } i<j \\
C_{i-j,1} C_{i-j+1,2} \dotsb C_{i-1,j-1}C_{i,j} & \text{ if } i\geq j.
\end{array}
\end{equation}
Furthermore, all cpo's are of this form, hence there are
$\frac{n(n+1)}{2}$ of them.
\end{lemma}

\begin{proof}
First note that since $i+j\leq n$, the support of $C_{i,j}$, namely
$I_{i,j}= \{s \in \widetilde{S} : s\leq C_{i,j}\}$, is a proper subset
of $\tilde{S}$, and the support of any segment smaller than $C_{i,j}$
in weak order is contained in $I_{i,j}$.  Any allowed product of
segments with support in $I_{i,j}$ would have length no longer than
the products in \eqref{e:cpo.type} by Lemma~\ref{l:allowed.a}.
Therefore, by Theorem~\ref{t:surjective.An}, the product above is the
unique cpo with support $I_{i,j}$.

Conversely, given an arbitrary cpo $w$ with support $I \subset
\tilde{S}$, $I \neq \tilde{S}$, let $r(w) = C_{i_{1},j_{1}}\dotsb
C_{i_p, j_p}$ be the corresponding reduced factorization into
segments.  If $\ell(C_{i_p, j_p}) >n$ then its support is all of
$\tilde{S}$ by definition of $C_{i_p, j_p}$, contradiction.  
Therefore, $i_{p}+j_{p}\leq n$ and $r(w)$ must be one of the products above. 
\end{proof}

\begin{lemma}\label{l:thin.not.pal.a}
Assume $W$ has type $A_{n}$ for $n\geq 2$.  

\begin{enumerate}
\item For each segment $C_{i,j}$ such that $i+j\leq n$ and each $k\geq 1$, the following elements are thin but not palindromic:
\begin{equation*}\label{e:thin.types}
\begin{array}{ll}
C_{0,j-i}C_{1,j-i+1}\dotsb C_{i-1,j-1}\cdot C_{i,j} \cdot  (C_{i+1,n-j})^{k} &  \text{if  }1\leq  i<j< n, \\
C_{0,j-i}C_{1,j-i+1}\dotsb C_{i-1,j-1}\cdot C_{i,j} \cdot (C_{n+1-i,j})^{k}&  \text{if  } 1\leq  i<j < n,\\
C_{i-j,1} C_{i+1,2} \dotsb C_{i-1,j-1}\cdot C_{i,j} \cdot (C_{i+1,n-j})^{k} &  \text{if  } 2\leq j\leq i < n, \\
C_{i-j,1} C_{i+1,2} \dotsb C_{i-1,j-1} \cdot C_{i,j} \cdot (C_{n+1-i,j})^{k} & \text{if  } 2\leq j\leq  i < n.
\end{array}
\end{equation*}
\item For each $0\leq i < n$, $1\leq j<n$, and $k>0$ such that $i+k+1
\neq 0 \ \mathrm{mod}\ n$ and $j+k \neq 0 \ \mathrm{mod}\ n$, the
elements $C_{i,1}(C_{n,1})^{k}$ and $C_{0,j}(C_{1,n})^{k}$ are thin
but not palindromic.

\item For each pair $1\leq i,j\leq n$ such that $i+j=n+1$, $k\geq 1$,
$k \neq 0\ \mathrm{mod}\ n$, the elements $(C_{i,j})^{k} $ are thin
but not palindromic.
\end{enumerate}
\end{lemma}

\begin{proof}
Using the relations on generators and segments from
Lemma~\ref{l:relations.An} and arguments similar to ones used in other
types, one can show that every element $w$ in the list above either
covers 2 or more elements in the Bruhat order or $w$ covers one
element $v$ and $v$ covers 3 elements.  Comparing with Bott's formula
\eqref{e:bott.formula} shows $w$ cannot be palindromic.

For example, consider elements of the form $w=(C_{1,n})^{k}$ for all
$k\geq 1$ and $k \neq 0\ \mathrm{mod}\ n$. First, assume that $1\leq
k<n$.  We claim $w$ covers $v=C_{1,n-k}(C_{2,n})^{k-1}$ and
$v'=C_{0,n}(C_{1,n})^{k-1}$.  To see the first covering relation,
recall $C_{1,n}=s_{1}s_{2}\dotsb s_{n}s_{0}$ and knock out the
leftmost $s_{k+1}$ generator from the left in $r(w)$.  The leftmost
$s_{k}$ then commutes to the right decreasing it's index as it passes
each $C_{1,n}$ until it becomes an $s_{2}$ which can glue onto the
rightmost $C_{1,n}$ forming a $C_{2,n}$.  Then the leftmost $s_{k-1}$
commutes right until it becomes an $s_{2}$ etc. Second, for the case
$k>n$, we make the following observation from the relations on
generators:
\[
C_{2,n}(C_{1,n})^{n} =(C_{1,n})^{n}   C_{2,n}
\]
and 
\[
C_{2,n}(C_{1,n})^{k} \not \in \minreps \text{ for } 1\leq k<n.
\]
Therefore, for any $w=(C_{1,n})^{k}$ with $k> n$ and $k \neq 0 \
\mathrm{mod}\ n$, we can apply the procedure in the first case to the
leftmost $(k\ \mathrm{mod}\ n)$ segments and commute the $C_{2,n}$'s
to the right past any factors of the form $(C_{1,n})^{n}$.  Hence, $w$
must actually cover at least 2 elements in Bruhat order, so it is not
palindromic.
\end{proof}

Let $\tchoose{n}{k}$ denote
the $t$-analog of $\chs{n}{k}$
\[
\tchoose{n}{k} = \frac{[n]!}{[k]! \ [n-k]!}
\]
where $[n]! = \prod_{i=0}^{n-1}(1+t+t^{2}+\dotsb + t^{i})$ is the
usual $t$-analog of $n!$.  The polynomial $\tchoose{n}{k}$ gives the
rank generating function for the partitions that fit inside a $k
\times n-k$ rectangle \cite[Sect. 1.3]{ec1}, hence it is palindromic.
The polynomial is also the Poincar\'e polynomial for the Grassmannian
manifold $G_{n,k}$.

The following lemma was originally proved in \cite{Mitchell-86}.  This
short proof was suggested by Hugh Thomas.  

\begin{lemma}\label{l:poincare.spiral}
Let $w$ be a spiral element of length $k\cdot n$ for $k\geq 1$.  
Then 
\[
P_{w}(t)= \tchoose{n+k}{k}.
\]
Hence $w$ is palindromic.
\end{lemma}

\begin{proof}
Let $w$ be the spiral element $C_{n-k-1,1}(C_{n,1})^{k}$ for some
$k\geq 1$. Then $\boundedtominreps^{-1}(w) = (n)^{k}$.  Every
partition $\mu \subset (n)^{k}$ has the property that
$\boundedtominreps(\mu ) \leq w$ in Bruhat order since the reduced
expression corresponding to the $n$-bounded partition $\mu$ is a
subexpression of $C_{n-k-1,1}(C_{n,1})^{k}$.  Therefore, to prove the
lemma, we only need to show that no $\lambda \in B^{n}$ exists such
that $\lambda \not \subset (n)^{k}$ and
$\boundedtominreps(\lambda)<w$.   

Given any $n$-bounded partition $\lambda$ such that $\lambda \not
\subset (n)^{k}$, then $\lambda $ must have more than $k$ parts.  If
the canonical reduced expression $\boundedtominreps(\lambda)=s_{i_{1}}\dotsb
s_{i_{p}}$, then the core partition $r_{i_{1}}\dotsb r_{i_{p}}
(\emptyset) = \bounded^{-1}(\lambda)$ must have at least $k+1$ rows.
However, by similar logic $\bounded^{-1}((n)^{k})$ has exactly $k$
rows, therefore $\boundedtominreps(\lambda) \not \leq
\boundedtominreps((n)^{k})=w$ by Lascoux's Theorem.  

The proof for $w=C_{0,n-k}(C_{1,n})^{k}$ follows from the natural
involution on the Dynkin diagram fixing $s_{0}$.  
\end{proof}

\begin{theorem}\label{t:palindromy.a}
Let $W$ be the Weyl group of type $A_{n}$.  Every element in
$\minreps$ is palindromic if $n=1$.  If $n\geq 2$, then the
palindromic elements in $\minreps$ are precisely the cpo's and the
spiral elements.
\end{theorem}

\begin{proof}
Lemma~\ref{l:cpos.a} and Lemma~\ref{l:poincare.spiral} prove each of
the cpo's and spiral elements are palindromic.  Lemmas~\ref{l:thin.a}
and~\ref{l:thin.not.pal.a} show that every other element is not
palindromic.
\end{proof}

\begin{corollary}\label{c:mw.a}
For $W$ of type $A_{n}$ with $n\geq 2$, we have $m_{W}=2$.  
\end{corollary}

\section{Appendix}\label{s:appendix}

\noindent \textbf{Example in type $F_{4}$}: The segments for $F_{4}$
are as follows:

\begin{equation}\label{e:segments.f4}
\begin{array}{lc}
\text{segment } & \text{length }
\\
\hrulefill & \hrulefill \\
 \segment{}{}(1) = s_0 & 1 \\
 \segment{}{}(2) = s_1s_0 & 2 \\
 \segment{}{}(3) = s_2s_1s_0 & 3 \\
 \segment{}{}(4) = s_3s_2s_1s_0 & 4 \\
 \segment{}{}(5) = s_2s_3s_2s_1s_0 & 5 \\
 \segment{1}{}(5) = s_4s_3s_2s_1s_0 & 5 \\
 \segment{}{}(6) = s_1s_2s_3s_2s_1s_0 & 6 \\
 \segment{1}{}(6) = s_2s_4s_3s_2s_1s_0 & 6 \\
 \segment{}{}(7) = s_1s_2s_4s_3s_2s_1s_0 & 7 \\
 \segment{1}{}(7) = s_3s_2s_4s_3s_2s_1s_0 & 7 \\
 \segment{}{}(8) = s_1s_3s_2s_4s_3s_2s_1s_0 & 8 \\
 \segment{1}{}(8) = s_2s_3s_2s_4s_3s_2s_1s_0 & 8 \\
 \segment{}{}(9) = s_1s_2s_3s_2s_4s_3s_2s_1s_0 & 9 \\
 \segment{1}{}(9) = s_2s_1s_3s_2s_4s_3s_2s_1s_0 & 9 \\
 \segment{}{}(10) = s_2s_1s_2s_3s_2s_4s_3s_2s_1s_0 & 10 \\
 \segment{1}{}(10) = s_3s_2s_1s_3s_2s_4s_3s_2s_1s_0 & 10 \\
 \segment{}{}(11) = s_3s_2s_1s_2s_3s_2s_4s_3s_2s_1s_0 & 11 \\
 \segment{1}{}(11) = s_4s_3s_2s_1s_3s_2s_4s_3s_2s_1s_0 & 11 \\
 \segment{}{}(12) = s_2s_3s_2s_1s_2s_3s_2s_4s_3s_2s_1s_0 & 12 \\
 \segment{1}{}(12) = s_4s_3s_2s_1s_2s_3s_2s_4s_3s_2s_1s_0 & 12 \\
 \segment{}{}(13) = s_2s_4s_3s_2s_1s_2s_3s_2s_4s_3s_2s_1s_0 & 13 \\
 \segment{}{}(14) = s_3s_2s_4s_3s_2s_1s_2s_3s_2s_4s_3s_2s_1s_0 & 14 \\
 \segment{}{}(15) = s_2s_3s_2s_4s_3s_2s_1s_2s_3s_2s_4s_3s_2s_1s_0 & 15 \\
 \segment{}{}(16) = s_1s_2s_3s_2s_4s_3s_2s_1s_2s_3s_2s_4s_3s_2s_1s_0 & 16 \\
\end{array}
\end{equation}

Observe that the segments do not have distinct lengths, so colors will
be necessary for the corresponding partitions.   

One can verify that the product of two segments $\segment{}{}(i)\cdot
\segment{}{}(j)$ is a reduced expression for a minimal length coset
representative in $\minreps$ if and only if the ordered pair $(i.\ j)$
appears in the following list of \textit{allowed patterns}:

\begin{small}
\[
\begin{array}{l}
 (16.16)
\\
 (15.16) (15.15)
\\
 (14.16) (14.15) (14.14)
\\
 (13.16) (13.15) (13.14)
\\
 (12.16) (12.15) (12.14)
\\
 (12^1.16) (12^1.15) (12^1.14) (12^1.13)
\\
 (11.16) (11.15) (11.14) (11.13)
\\
 (11^1.16) (11^1.15) (11^1.14) (11^1.13) (11^1.12^1) (11^1.11^1)
\\
 (10.16) (10.15) (10.14) (10.13)
\\
 (10^1.16) (10^1.15) (10^1.14) (10^1.13) (10^1.12^1) (10^1.11^1)
\\
 (9.16) (9.15) (9.14) (9.13) (9.12)
\\
 (9^1.16) (9^1.15) (9^1.14) (9^1.13) (9^1.12^1) (9^1.11^1)
\\
 (8.16) (8.15) (8.14) (8.13) (8.12^1) (8.11^1)
\\
 (8^1.16) (8^1.15) (8^1.14) (8^1.13) (8^1.12)
\\
 (7.16) (7.15) (7.14) (7.13) (7.12^1) (7.11^1)
\\
 (7^1.16) (7^1.15) (7^1.14) (7^1.13) (7^1.12) (7^1.12^1) (7^1.11^1)
\\
 (6.16) (6.15) (6.14) (6.13) (6.12^1) (6.11^1)
\\
 (6^1.16) (6^1.15) (6^1.14) (6^1.13) (6^1.12) (6^1.12^1) (6^1.11^1)
\\
 (5.16) (5.15) (5.14) (5.13) (5.12) (5.12^1) (5.11^1)
\\
 (5^1.16) (5^1.15) (5^1.14) (5^1.13) (5^1.12) (5^1.12^1) (5^1.11) (5^1.11^1) (5^1.10^1)
\\
 (4.16) (4.15) (4.14) (4.13) (4.12) (4.12^1) (4.11) (4.11^1) (4.10^1)
\\
 (3.16) (3.15) (3.14) (3.13) (3.12) (3.12^1) (3.11) (3.11^1) (3.10^1)
\\
 (2.16) (2.15) (2.14) (2.13) (2.12) (2.12^1) (2.11) (2.11^1) (2.10) (2.10^1) (2.9^1)
\\
 (1.16) (1.15) (1.14) (1.13) (1.12) (1.12^1) (1.11) (1.11^1) (1.10) (1.10^1) (1.9) (1.9^1) (1.8) (1.7) (1.6)
\end{array}
\]
\end{small}

This list of allowed patterns is harder to internalize than in types
$B,C,D,G_{2}$, however some useful information pops out of it.  First,
the only parts which can repeat in a colored partition are
$11^{1},14,15$ and $16$.  The parts indicated by $14,15,16$ can occur
with all other parts, but $11^{1}$ cannot appear in any partition with
$8^{1},9,10,11$ or $12$.  Let $NR$ be the set of allowed partitions
with no parts in the set $\{11^{1},14,15,16 \}$.  There are 132
partitions in $NR$ and they all have size at most 40.  The colored
partition $(13,12^{1},10^{1},5^{1})$ is the unique one in $NR$ with size 40
exactly.  Let $NRE$ be the set of all partitions in $NR$ that can be
in an allowed partition with $11^{1}$.  There are 96 elements in $NRE$
including $(13,12^{1},10^{1},5^{1})$.  

To compute the generating function for $\partitions(F_{4})$, let
\[
G_{NR}(t) = \sum_{\lambda \in NR} t^{|\lambda|} 
\]
be the rank generating function for the partitions in $NR$ and similarly for $NRE$.
Then, since every partition in $\partitions(F_{4})$ contains an $11^{1}$ or not, we get the following identity
\begin{equation}\label{e:identify.f4}
G_{\partitions(F_{4})} = \frac{G_{NR}}{(1-t^{14})(1-t^{15})(1-t^{16})}
+ \frac{t^{11} \cdot
G_{NRE}}{(1-t^{11})(1-t^{14})(1-t^{15})(1-t^{16})}.  
\end{equation}
By computer it is easy to compute the coefficients for $G_{NR}$ and $G_{NRE}$ 
\[
\begin{array}{ll}
NR &	(1\ 1\ 1\ 1\ 1\ 2\ 2\ 3\ 3\ 3\ 4\ 4\ 5\ 5\ 5\ 5\ 4\ 5\ 5\ 6\ 6\ 6\ 6\ 6\ 5\ 5\ 4\ 4\ 3\ 2\ 2\ 2\ 3\ 2\ 2
\ 2\ 2\ 1\ 1\ 1\ 1)\\
NRE &	(1\ 1\ 1\ 1\ 1\ 2\ 2\ 3\ 2\ 2\ 2\ 2\ 2\ 3\ 3\ 3\ 2\ 3\ 4\ 5\ 5\ 4\ 4\ 4\ 3\ 3\ 3\ 3\ 2\ 1\ 2\ 2\ 3\ 2\ 2
\ 2\ 2\ 1\ 1\ 1\ 1)
\end{array}
\]
Miraculously using Maple, the sum in \eqref{e:identify.f4} simplifies
to
\[
\frac{1}{(1-t^{1})(1-t^{5})(1-t^{7}) (1-t^{11})}.
\]
as expected using Bott's formula since the exponents for $F_{4}$ are
${1,5,7,11}$.  This proves Theorem~\ref{t:partition.bijection} for
type $F_{4}$.


The nine palindromic elements of $F_{4}$ are given by the affine partitions
\[
(),(1),(2),(3),(4),(5),(5^{1}), (6), (6,1).
\]
Using \eqref{e:segments.f4}, we see the cpo's correspond with
$(),(1),(2),(3), (6,1).$

\newpage

\newcommand{\littlespace}{.2in}
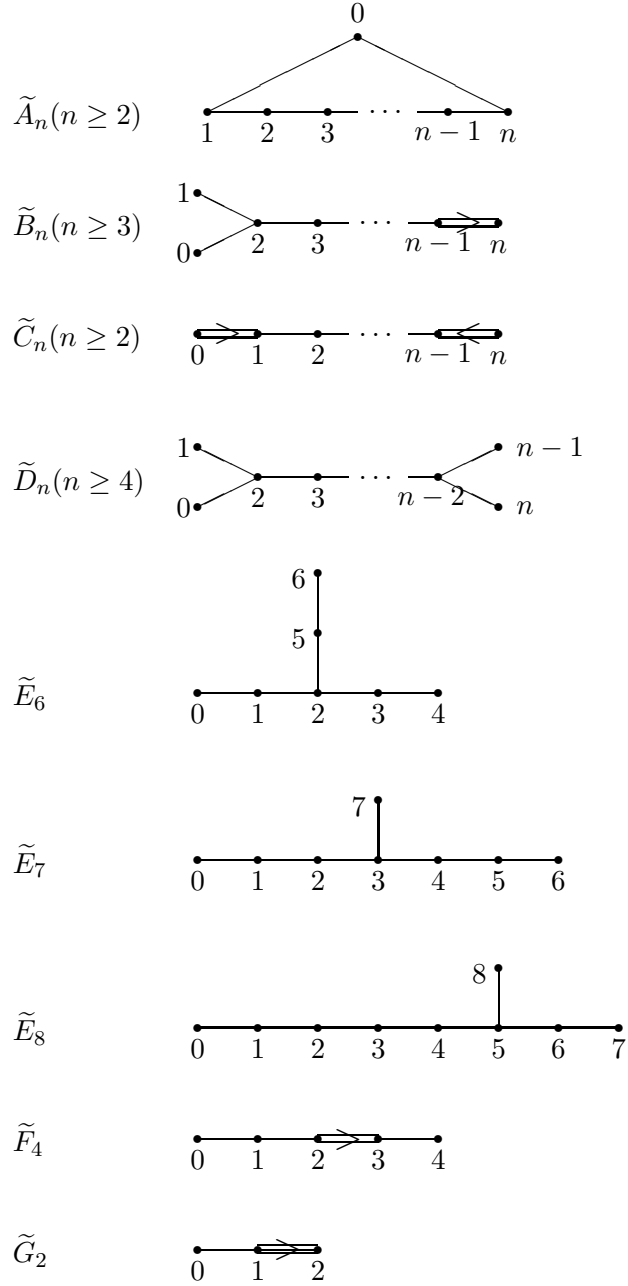
\begin{figure}\label{f:dynkin}
\caption{Affine Dynkin diagrams}
\setlength{\unitlength}{.8cm}
$$
\begin{array}{lll}
\\ \\ \vspace{\littlespace}
\widetilde{A}_n (n\geq 2)&& 
\raisebox{1ex}{
\begin{picture}(5,0)
\mp(0,0)(1,0){3}{\ci}
\put(0,0){\num{1}}\put(1,0){\num{2}}\put(2,0){\num{3}}
\put(3,0){\makebox(0,0){\ldots}}
\mp(5,0)(-1,0){2}{\ci}
\put(4,0){\num{n-1}}\put(5,0){\num{n}}
\put(0,0){\line(1,0){2.5}}
\put(5,0){\line(-1,0){1.5}}
\put(2.5,2){\num{0}}
\put(0,0){\line(2,1){2.5}}
\put(5,0){\line(-2,1){2.5}}
\put(2.5,1.25){\ci}
\end{picture}
}
\\ \\ \vspace{\littlespace}
\widetilde{B}_n (n\geq 3)&& 
\raisebox{1ex}{\begin{picture}(5,0)
\mp(1,0)(1,0){2}{\ci}
\put(1,0){\num{2}}\put(2,0){\num{3}}
\put(3,0){\makebox(0,0){\ldots}}
\mp(5,0)(-1,0){2}{\ci}
\put(4,0){\num{n-1}}\put(5,0){\num{n}}
\put(1,0){\line(1,0){1.5}}
\put(4,0){\line(-1,0){.5}}
\put(4,.06){\line(1,0){1}}
\put(4,-.06){\line(1,0){1}}
\put(4.5,0){\makebox(0,0){\Large$>$}}
\put(0,.5){\ci}\put(0,-.5){\ci}
\put(0,.5){\makebox(0,0)[l]{\hspace{-.35\unitlength}$1$}}
\put(0,-.5){\makebox(0,0)[l]{\hspace{-.35\unitlength}$0$}}
\put(1,0){\line(-2,1){1}}
\put(1,0){\line(-2,-1){1}}
\end{picture}}
\\ \\ \vspace{\littlespace}
\widetilde{C}_n (n\geq 2)&& 
\raisebox{1ex}{\begin{picture}(5,0)
\mp(0,0)(1,0){3}{\ci}
\put(1,0){\num{1}}\put(2,0){\num{2}}
\put(3,0){\makebox(0,0){\ldots}}
\mp(5,0)(-1,0){2}{\ci}
\put(4,0){\num{n-1}}\put(5,0){\num{n}}
\put(1,0){\line(1,0){1.5}}
\put(4,0){\line(-1,0){.5}}
\put(4,.06){\line(1,0){1}}
\put(4,-.06){\line(1,0){1}}
\put(4.5,0){\makebox(0,0){\Large$<$}}
\put(.5,0){\makebox(0,0){\Large$>$}}
\put(0,.06){\line(1,0){1}}
\put(0,-.06){\line(1,0){1}}
\put(0,0){\num{0}}
\end{picture}}
\\ \\ \vspace{\littlespace}
\widetilde{D}_n (n\geq 4)&& 
\raisebox{1ex}{\begin{picture}(5,.8)
\mp(1,0)(1,0){2}{\ci}
\put(1,0){\num{2}}\put(2,0){\num{3}}
\put(3,0){\makebox(0,0){\ldots}}
\put(4,0){\ci}
\put(4,0){\num{n-2\,\,\,}}
\put(1,0){\line(1,0){1.5}}
\put(4,0){\line(-1,0){.5}}
\put(4,0){\line(2,1){1}}
\put(4,0){\line(2,-1){1}}
\put(5,.5){\ci}\put(5,-.5){\ci}
\put(5,.5){\makebox(0,0)[l]{\hspace{.3\unitlength}$n-1$}}
\put(5,-.5){\makebox(0,0)[l]{\hspace{.3\unitlength}$n$}}
\put(0,.5){\ci}\put(0,-.5){\ci}
\put(0,.5){\makebox(0,0)[l]{\hspace{-.35\unitlength}$1$}}
\put(0,-.5){\makebox(0,0)[l]{\hspace{-.35\unitlength}$0$}}
\put(1,0){\line(-2,1){1}}
\put(1,0){\line(-2,-1){1}}
\end{picture}}
\\ \\ \vspace{\littlespace}
\widetilde{E}_{6} && 
\raisebox{1ex}{\begin{picture}(5,2)
\mp(0,0)(1,0){5}{\ci}
\mp(2,1)(0,1){2}{\ci}
\put(0,0){\num{0}}\put(1,0){\num{1}}\put(2,0){\num{2}}\put(3,0){\num{3}}
\put(4,0){\num{4}}
\put(2,.9){\makebox(0,0)[r]{$5$\hspace{.2\unitlength}}}
\put(2,1.9){\makebox(0,0)[r]{$6$\hspace{.2\unitlength}}}
\put(2,0){\line(0,1){2}}
\put(0,0){\line(1,0){4}}
\end{picture}}
\\ \\ \vspace{\littlespace}
\widetilde{E}_{7} && 
\raisebox{1ex}{\begin{picture}(5,1.2)
\mp(0,0)(1,0){7}{\ci}\put(3,1){\ci}
\put(0,0){\num{0}}\put(1,0){\num{1}}\put(2,0){\num{2}}\put(3,0){\num{3}}
\put(4,0){\num{4}}\put(5,0){\num{5}}\put(6,0){\num{6}}
\put(3,.9){\makebox(0,0)[r]{$7$\hspace{.2\unitlength}}}
\put(3,0){\line(0,1){1}}
\put(0,0){\line(1,0){6}}

\end{picture}}
\\ \\ \vspace{\littlespace}
\widetilde{E}_{8} && 
\raisebox{1ex}{\begin{picture}(5,1.2)
\mp(0,0)(1,0){8}{\ci}\put(5,1){\ci}
\put(0,0){\num{0}}\put(1,0){\num{1}}\put(2,0){\num{2}}\put(3,0){\num{3}}
\put(4,0){\num{4}}\put(5,0){\num{5}}\put(6,0){\num{6}} \put(7,0){\num{7}} 
\put(5,.9){\makebox(0,0)[r]{$8$\hspace{.2\unitlength}}}
\put(5,0){\line(0,1){1}}
\put(0,0){\line(1,0){7}}
\end{picture}}
\\ \\ \vspace{\littlespace}
\widetilde{F}_4 && 
\raisebox{1ex}{\begin{picture}(5,0)
\mp(0,0)(1,0){5}{\ci}
\put(0,0){\num{0}}
\put(1,0){\num{1}}\put(2,0){\num{2}}\put(3,0){\num{3}}\put(4,0){\num{4}}
\put(0,0){\line(1,0){2}}
\put(3,0){\line(1,0){1}}
\put(2,.06){\line(1,0){1}}
\put(2,-.06){\line(1,0){1}}
\put(2.5,0){\makebox(0,0){\Large$>$}}
\end{picture}}
\\ \\ \vspace{\littlespace}
\widetilde{G}_2 && 
\raisebox{1ex}{\begin{picture}(5,0)
\mp(0,0)(1,0){3}{\ci}
\put(0,0){\num{0}}\put(1,0){\num{1}}\put(2,0){\num{2}}
\put(0,0){\line(1,0){2}}
\put(1,.08){\line(1,0){1}}
\put(1,-.06){\line(1,0){1}}
\put(1.5,0.01){\makebox(0,0){\Large$>$}}
\end{picture}}
\\
\end{array}
$$
\end{figure}

\newpage

\section{Acknowledgments} Many thanks to Andrew Crites, Adriano
Garsia, Thomas Lam, Mark Shimozono for helpful suggestions to this
work.  Special thanks to Hugh Thomas for the short proof of
Lemma~\ref{l:poincare.spiral}.

\bibliographystyle{siam}

\def\cprime{$'$}

\end{document}